\documentclass[12pt]{amsart}

\usepackage{color,graphicx}
\usepackage{psfrag}
\usepackage{amssymb}                    
\usepackage{amsmath}                    
\usepackage{amsthm}                     
\usepackage{latexsym}                   
\usepackage{subfigure}
\usepackage{epsfig}                     
\usepackage{pstricks,pst-node,pst-text,pst-tree}
\usepackage{amscd}
\usepackage{amsmath}
\usepackage{psfrag}
\usepackage{xypic}

\newtheorem{theorem}{Theorem}[section]
\newtheorem{corollary}[theorem]{Corollary}
\newtheorem{lemma}[theorem]{Lemma}
\newtheorem{conjecture}[theorem]{Conjecture}
\newtheorem{proposition}[theorem]{Proposition}

\theoremstyle{definition}
\newtheorem{definition}[theorem]{Definition}

 \newtheorem{example}[theorem]{Example}
 \newtheorem{remark}[theorem]{Remark}
\newtheorem{algorithm}[theorem]{Algorithm}

\DeclareMathOperator{\Spec}{Spec}
\DeclareMathOperator{\Proj}{Proj}

\DeclareMathOperator{\face}{face}

\newcommand{\one}{\ensuremath{(\mathrm{i})}}
\newcommand{\two}{\ensuremath{(\mathrm{ii})}}

\newcommand{\C}{\ensuremath{\mathbb{C}}}

 \newcommand{\GL}{\ensuremath{\operatorname{GL}}}
 \newcommand{\Hom}{\ensuremath{\operatorname{Hom}}}
 
 \newcommand{\K}{\Bbbk}

 \newcommand{\N}{\ensuremath{\mathbb{N}}}

 \newcommand{\Q}{\ensuremath{\mathbb{Q}}}
  
 \newcommand{\SL}{\ensuremath{\operatorname{SL}}}

 \newcommand{\Z}{\ensuremath{\mathbb{Z}}}

 \newcommand{\conv}{\ensuremath{\operatorname{conv}}}
 \newcommand{\diag}{\ensuremath{\operatorname{diag}}}
 
 \newcommand{\ghilb}{\ensuremath{G}\operatorname{-Hilb}}
 \newcommand{\git}{\ensuremath{\operatorname{/\!\!/}}}
 \newcommand{\hilb}{\operatorname{Hilb}}
  \newcommand{\hilbg}{\operatorname{Hilb}^{\ensuremath{G}}}
 
 \newcommand{\owe}{\ensuremath{\mathcal{O}}}

 \newcommand{\st}{\ensuremath{\operatorname{\bigm{|}}}}

 \numberwithin{equation}{section}
\begin{document}

 \bibliographystyle{plain} 
 
 \title[Coherent component for McKay quiver representations]{Moduli of McKay quiver representations I: \\ the coherent component}


 %
 \author{Alastair Craw} \address{Department of Mathematics, University of Glasgow, Glasgow G12 8QW, UK.} \email{craw@maths.gla.ac.uk}
 
 \author{Diane Maclagan} \address{Department of Mathematics,
       Hill Center-Busch Campus,
       Rutgers University,
       110 Frelinghuysen Rd,
       Piscataway, NJ 08854, USA.} \email{maclagan@math.rutgers.edu}

 \author{Rekha R.\ Thomas} \address{Department of Mathematics, University
   of Washington, Seattle, WA
   98195, USA.}\email{thomas@math.washington.edu}


 \begin{abstract}
   For a finite abelian group \(G\subset \GL(n,\K)\), we describe the
   coherent component \(Y_\theta\) of the moduli space
   \(\mathcal{M}_\theta\) of $\theta$-stable McKay quiver
   representations.  This is a not-necessarily-normal toric variety
   that admits a projective birational morphism \(Y_\theta \rightarrow
   \mathbb{A}_{\K}^n/G\) obtained by variation of GIT quotient.  As a
   special case, this gives a new construction of Nakamura's
   $G$-Hilbert scheme $\hilbg$ that avoids the (typically highly
   singular) Hilbert scheme of $\vert G\vert$-points in
   $\mathbb{A}^n_{\K}$.  To conclude, we describe the toric fan of
   \(Y_\theta\) and hence calculate the quiver representation
   corresponding to any point of \(Y_\theta\).
 \end{abstract}

 \maketitle

 \section{Introduction}
 In this paper we concretely describe the coherent component
 \(Y_\theta\) of the moduli spaces \(\mathcal{M}_\theta\) of
 representations of the McKay quiver for a finite abelian subgroup
 \(G\subset \GL(n,\K)\) and generic parameter $\theta$, where $\K$ is
 an algebraically closed field whose characteristic does not divide
 $r:= \vert G\vert$.  The irreducible component $Y_\theta$ is a
 not-necessarily-normal toric variety that admits a projective
 birational morphism \(Y_\theta \rightarrow \mathbb{A}_{\K}^n/G\)
 obtained by variation of GIT quotient.

 The motivation to study the moduli spaces \(\mathcal{M}_\theta\)
 comes from their role in the McKay correspondence (see \cite{BKR,
 Haiman, Reid3, BezrukavnikovKaledin}).  For a finite subgroup
 \(G\subset \SL(n,\K)\), this is the expected equivalence between the
 $G$-equivariant geometry of \(\mathbb{A}^n_{\K}\) and the geometry of
 a crepant resolution \(Y\rightarrow \mathbb{A}^n_{\K}/G\) (if one
 exists).  For \(n\leq 3\), Kronheimer~\cite{Kronheimer} and
 Bridgeland--King--Reid~\cite{BKR} proved that the moduli spaces
 \(\mathcal{M}_\theta\) are crepant resolutions of
 $\mathbb{A}^n_{\K}/G$ for all generic parameters \(\theta\in
 \Theta\), where $\Theta$ is a rational vector space of weights.
 This moduli interpretation of the crepant resolution enabled
 \cite{BKR} to establish the McKay correspondence as an equivalence of
 derived categories for \(n\leq 3\).  Craw--Ishii \cite{CrawIshii}
 established a partial converse for finite abelian subgroups $G
 \subset \SL(3,\K)$: every projective crepant resolution of
 $\mathbb{A}^3_{\K}/G$ is isomorphic to $\mathcal M_{\theta}$ for some
 generic parameter $\theta\in \Theta$.  For $n\geq 4$, it is unknown
 whether every projective crepant resolution of
 \(\mathbb{A}^n_{\K}/G\) can be constructed as (a component of)
 $\mathcal M_{\theta}$ for some generic \(\theta\in \Theta\). The
 tools introduced in this paper allow the investigation of such
 questions.  For example, we show in Example~\ref{ex:weightoneaction}
 below that, for the diagonal action of the group $G = \Z/n\Z$ on
 $\mathbb{A}^n_{\K}$ with weights $(1,\dots,1)$, the unique toric
 crepant resolution $Y$ of $\mathbb{A}^n_{\K}/G$ is isomorphic to the
 component $Y_\theta$ of $\mathcal M_{\theta}$ for any generic
 \(\theta\in \Theta\).

 Since $G$ is abelian, we may assume that $G$ is contained in the
 subgroup $(\K^*)^n$ of diagonal matrices with nonzero entries in
 $\GL(n,\K)$.  This representation of $G$ decomposes into irreducible
 representations \(\rho_1\oplus \dots \oplus\rho_n\).  The \emph{McKay
   quiver} is the directed graph whose vertices are the irreducible
 representations $\rho$ of $G$, with an arrow from $\rho\rho_i$ to
 $\rho$ for every $\rho$ and \(1\leq i\leq n\).  Since \(G\) is
 abelian, there are \(r\) vertices and \(nr\) arrows.  We consider
 McKay quiver representations of dimension vector $(1,\dots,1)$, which
 correspond to points in $\mathbb{A}^{nr}_{\K}$.  Requiring certain
 commutativity relations gives a subscheme \(Z\subset
 \mathbb{A}^{nr}_{\K}\).  An algebraic torus \(T_B = (\K^*)^r/\K^*\)
 acts by change of basis on quiver representations, which in turn
 gives an action of \(T_B\) on \(Z\).  Moduli spaces of
 representations satisfying the commutativity relations are
 constructed by Geometric Invariant Theory as quotients \(\mathcal
 M_{\theta} := Z\git_{\theta}T_B\), where \(\theta\in \Theta \cong
 \Q^{r-1}\) is a fractional character of $T_B$.  The best known
 example is the $G$-Hilbert scheme \(\ghilb\) parametrizing
 \emph{\(G\)-clusters} (see~\cite{ItoNakajima}).

 The first main result of this paper constructs explicitly the
 component \(Y_\theta\) of \(\mathcal{M}_\theta\) that is birational
 to \(\mathbb{A}^n_{\K}/G\).  The crucial step is to introduce an
 \((r+n)\times nr\)-matrix $C$ obtained by augmenting the vertex-edge
 incidence matrix of the McKay quiver. This matrix defines an
 irreducible component $V$ of the scheme $Z\subset
 \mathbb{A}^{nr}_{\K}$ that has the following properties (see
 Theorems~\ref{t:cohcompeqtns} and~\ref{thm:Ytheta}).
 
 \begin{theorem} 
 \label{thm:1.1}
 The scheme $Z$ has a unique irreducible component $V = \Spec \K[\mathbb N C]$
 that does not lie in any coordinate hyperplane of
 $\mathbb{A}^{nr}_{\K}$, where $\mathbb N C$ is the semigroup
 generated by the columns of the matrix $C$.  In addition:
 \begin{enumerate}
 \item For $\theta \in \Theta$, the GIT quotient $Y_\theta :=
   V\git_\theta T_B$ is a not-necessarily-normal toric variety that
   admits a projective birational morphism $\tau_\theta\colon Y_\theta
   \rightarrow \mathbb{A}^{n}_{\K}/G$ obtained by variation of GIT
   quotient.
 \item For generic $\theta\in \Theta$, the variety $Y_\theta$ is the
   unique irreducible component of $\mathcal M_{\theta}$ containing
   the $T_B$-orbit closures of the points of $Z\cap (\K^*)^{nr}$.
 \end{enumerate}
 We call \(Y_{\theta}\) the \emph{coherent component} of the moduli
 space \(\mathcal M_{\theta}\).
 \end{theorem}

 \begin{corollary}
 \label{cor:1.2}
 For any \(\theta\in \Theta\) such that
   \(\mathcal{M}_\theta \cong \ghilb\), the coherent component
   \(Y_\theta\) is isomorphic to the irreducible scheme \(\hilbg\)
   introduced by Nakamura.  
 \end{corollary}

This corollary (see Proposition~\ref{prop:hilbg}) provides a direct
GIT construction of the scheme \(\hilbg\) that avoids the original
construction as a subscheme of the Hilbert scheme of \(r\)-points in
\(\mathbb{A}^n_{\K}\).  In addition, it confirms the suggestion of
Mukai that $\hilbg$ should be obtained from \(\mathbb{A}^n_{\K}/G\) by
variation of GIT quotient.

\medskip

The second main result describes explicitly the set of
$\theta$-semistable McKay quiver representations corresponding to
points of $Y_\theta$. For generic $\theta\in \Theta$, these
representations encode the restriction to $Y_\theta$ of the universal
quiver representation on the fine moduli space $\mathcal{M}_\theta$
(see Craw--Ishii~\cite[\S2]{CrawIshii}).  As a first step, we prove
that the toric fan of the variety $Y_{\theta}$ (see
Section~\ref{sec:toricfan} for details) is the inner normal fan of a
polyhedron \(P_\theta\) obtained by slicing the cone \(P\subseteq
\Q^{r+n}\) generated by the column vectors of the matrix \(C\) (see
Corollary~\ref{cor:fan-of-normalization}).  This explicit description
enables us to calculate new examples in detail.  In particular, we
obtain the new moduli description of the unique crepant toric
resolution $Y\rightarrow \mathbb{A}^n_{\K}/G$ for the action of $G =
\Z/n\Z$ on $\mathbb{A}^n_{\K}$ with weights $(1,\dots,1)$ mentioned
above.

Any vector ${\bf w}$ in the support of the toric fan defining
$Y_\theta$ determines a distinguished point of $Y_\theta$, and hence a
distinguished $\theta$-semistable McKay quiver representation
$b_{\theta,{\bf w}}\in\mathbb{A}^{nr}_{\K}$.  To compute the
coordinates of $b_{\theta,{\bf w}}$ explicitly in terms of $\theta$ and $\mathbf{w}$, consider the
slice \(P^\vee_{\bf w} := \{\mathbf{v} \in (\Q^{r})^* : w_i+v_\rho
- v_{\rho\rho_i}\geq 0\}\) of the polyhedral cone \(P^\vee\) dual to
\(P\).  The following result is proved in
Theorem~\ref{thm:distinguishedcone}.

 \begin{theorem}
 \label{thm:1.3}
 Fix $\theta$ and $\mathbf{w}$, and let $\mathbf{v} \in P^{\vee}_{\mathbf{w}}$ be any
 vector satisfying $\theta \cdot \mathbf{v} \leq \theta \cdot
 \mathbf{v}'$ for all $\mathbf{v}' \in P^{\vee}_{\mathbf{w}}$.  Then
 the distinguished $\theta$-semistable quiver representation
 $b_{\theta,{\bf w}}=(b_i^{\rho})$ has
 \begin{equation}
 b_i^\rho =
 \left\{\begin{array}{cl} 1 & \text{if } w_i+v_\rho - v_{\rho\rho_i} =
 0 \\ 0 & \text{if } w_i+v_\rho - v_{\rho\rho_i} > 0
 \end{array}\right. .
\end{equation} 

\end{theorem} 
Computing the representations $b_{\theta,{\bf w}}$ that give torus-fixed
points of $Y_\theta$ has been a key tool in understanding the moduli
spaces \(\mathcal{M}_\theta\) (see \cite{Reid2, Nakamura, Crawthesis,
  CrawIshii}), though no good algorithm was known in general until
now.

\medskip

From the perspective of string theory, the results of this paper are
as follows.  For \(\K = \C\), the spaces \(\mathcal{M}_\theta\) appear
in the physics literature as moduli of \(D0\)-branes on the orbifold
\(\C^n/G\), where \(\theta\) is a Fayet-Iliopoulos term for ${\rm
  U}(1)$ gauge multiplets present in the world-volume theory (see
\cite{DGM}).
 The matrix \(C\) introduced in Section~\ref{sec:cohcomp} encodes both
 the \(D\)-term equations and the \(F\)-term equations of the relevant
 quiver gauge theory.  More precisely, the top \(r\times nr\)
 submatrix \(B\) encodes the \(D\)-terms, giving the moment map for
 the action of ${\rm U}(1)^r/{\rm U}(1)$, and the bottom \(n\times
 nr\) submatrix of \(C\) encodes the \(F\)-terms obtained from the
 partial derivatives of the superpotential of the quiver gauge theory.
 
 We now explain the division into sections.  Section~\ref{sec:moduli}
 reviews the construction of the moduli spaces $\mathcal{M}_{\theta}$,
 including some well-known facts from Geometric Invariant Theory.
 Section~\ref{sec:cohcomp} introduces the irreducible component $V$ of
 $Z$.  Section~\ref{sec:toricGIT} constructs the coherent component
 $Y_{\theta}$ to complete the proof of Theorem~\ref{thm:1.1}, and
 Section~\ref{sec:hilbg} establishes Corollary~\ref{cor:1.2}.  The
 toric fan of $Y_\theta$ is computed in Section~\ref{sec:toricfan},
 and Theorem~\ref{thm:1.3} is established in Section~\ref{sec:distinguished}.

 \medskip

 \noindent \textbf{Conventions} For an integer matrix \(C\), let $\N
 C$ denote the semigroup generated by the columns of \(C\). Similarly,
 $\Z C$ denotes the lattice, $\Q_{\geq 0} C$ the rational cone and
 \(\Q C\) the rational vector space generated by columns of \(C\).
 A point on a scheme over \(\K\) means a closed point.  Write \(\K^*\)
 for the one-dimensional algebraic torus.

 \medskip

 \noindent \textbf{Acknowledgements} We thank Bernd
 Sturmfels for bringing us together.  The original observation of a
 link between $\ghilb$ and the toric Hilbert scheme is due to him.  We
 also thank Iain Gordon, Mark Haiman, Akira Ishii, S.~Paul Smith and Bal\'{a}zs
 Szendr\H{o}i for useful discussions.  Finally, we thank the
 organizers of PCMI 2004 for providing a stimulating environment where
 part of this paper was written.  The third author was partially
 supported by NSF grant DMS-04010147.

\medskip

\section{Moduli spaces of McKay quiver representations}
 \label{sec:moduli} In this section we recall the construction of the
 moduli spaces of representations of the McKay quiver for a finite
 abelian subgroup \(G\subset\GL(n,\K)\) of order \(r\), where \(\K\)
 is an algebraically closed field whose characteristic does not divide
 \(r\).  Since $G$ is abelian, we may assume that $G$ is contained in
 the subgroup $(\K^*)^n$ of diagonal matrices with nonzero entries in
 $\GL(n,\K)$.
 
Irreducible representations of $G$ are one-dimensional and hence
 define elements of the dual group of characters \(G^*\!
 :=\Hom(G,\K^*)\).  The \(n\)-dimensional representation given by the
 inclusion of \(G\) in \(\GL(n,\K)\) decomposes into one-dimensional
 representations \(\rho_1\oplus \dots \oplus \rho_n\) by Schur's
 lemma, so \(g\in G\) acts on \(\mathbb{A}_{\K}^n\) as the diagonal
 matrix \(\diag \!\big{(}\rho_1(g),\dots,\rho_n(g)\big{)}\).  Applying
 the functor $\Hom(-,\K^*)$ to the injective homomorphism of \(G\) into the
 algebraic torus \(T^n \cong (\K^{\ast})^n\) of diagonal matrices with
 nonzero entries induces a surjective homomorphism \begin{equation}
 \label{eqn:deg} \deg\colon \Z^n \rightarrow G^*, \end{equation} where
 \(\Z^n = \Hom(T^n,\K^*)\) is the character lattice of \(T^n\).  This
 determines a $G^*$-grading of the coordinate ring $\K[x_1,\dots,x_n]$
 of $\mathbb A^n_{\K}$ by setting
 $\deg(x_i):=\deg(\mathbf{e}_i)=\rho_i$, where $\mathbf{e}_i$ is a
 standard basis vector of $\mathbb Z^n$.  Since $\deg$ is surjective,
 $\rho_1,\dots, \rho_n$ generate the group $G^*$.

 \begin{definition}
 \label{defn:quiverrep}
 The \emph{McKay quiver} of $G\subset \GL(n,\K)$ is the directed graph with a
 vertex for each \(\rho\in G^*\), and an arrow \(a_i^\rho\) from
 \(\rho\rho_i\) to \(\rho\) for each \(\rho\in G^*\) and \(1\leq i
 \leq n\).  We say the arrow \(a_i^\rho\) is \emph{labeled} \(i\).
 \end{definition}

 This sign convention agrees with Sardo Infirri~\cite{SI2} but differs
 from that of Craw--Ishii~\cite{CrawIshii} where \(\rho_i\) is denoted
 \(\rho_i^{-1}\) and the arrows go from \(\rho\) to \(\rho\rho_i\).
 Note that the McKay quiver is strongly connected, since for any pair
 $\rho', \rho \in G^*$ there is a directed path from $\rho'$ to
 $\rho$.  Every such path comes from writing $\rho^{-1}\rho' =
 \bigotimes_{1\leq i\leq n} \rho_i^{m_i}$ for some $m_i \in \N$.

 \begin{example} \label{ex:runningexample} Consider the subgroup \(G\!
   := \Z/ 7\Z\subset \GL(2,\K)\) generated by the diagonal matrix \(g
   = \diag(\omega,\omega^2)\), where \(\omega\) is a primitive seventh
   root of unity.  This is the action of type \(\frac{1}{7}(1,2)\).
   Let $x:=x_1$, and $y:=x_2$ be the coordinates on $\mathbb
   A^2_{\K}$.  We have $\deg(x)=\rho_1$ and $\deg(y)=\rho_2$, where
   $\rho_1(g)=\omega$ and $\rho_2(g)=\omega^2$. The McKay quiver has
   vertices \(\sigma_0, \sigma_1, \dots, \sigma_6\) where $\sigma_0$
   is the trivial representation of $G$, and where $\sigma_1$ and
   $\sigma_2$ coincide with the representations $\rho_1$ and $\rho_2$
   respectively.  For each vertex $\sigma_j$ there are arrows
   $\sigma_{j+1}\rightarrow \sigma_j$ and $\sigma_{j+2}\rightarrow
   \sigma_j$ where addition is modulo 7. The arrows are denoted
   \(a^{\sigma_j}_i\) for \(i = 1,2\) and \(j = 0,\dots ,6\).  The
   quiver is shown in Figure~\ref{fig:mckayquiver}.

\begin{figure}[!ht]
\input{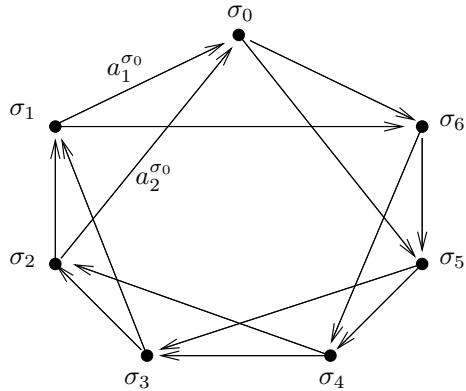}
 \caption{The McKay quiver for the action of type $\frac{1}{7}(1,2)$}
 \label{fig:mckayquiver}
\end{figure}
 \end{example}

 \begin{definition} \label{defn:quiverreps} A \emph{representation of
     the McKay quiver} with dimension vector \((1,\dots ,1)\in \N^r\)
   is the assignment of a one-dimensional \(\K\)-vector space
   \(R_\rho\) to each vertex \(\rho\), and a linear map
   \(R_{\rho\rho_i}\rightarrow R_{\rho}\) to each arrow $a_i^{\rho}$
   in the McKay quiver.  Fix a basis for each \(R_{\rho}\) and write
   \(b_i^{\rho} \in \K\) for the entry of the \(1 \times 1\) matrix of
   the linear map $R_{\rho\rho_i}\rightarrow R_{\rho}$.  We
   occasionally use $b_i^{\rho}$ to refer to the linear map itself.
 \end{definition}

\begin{figure}[!ht]
\psfrag{a}{$\rho \rho_i \rho_j$}
\psfrag{b}{$\rho \rho_j$}
\psfrag{c}{$\rho \rho_i$}
\psfrag{d}{$\rho$}
\psfrag{i}{$b_i^{\rho \rho_j}$}
\psfrag{k}{$b_j^{\rho \rho_i}$}
\psfrag{j}{$b_j^{\rho}$}
\psfrag{l}{$b_i^{\rho}$}
\epsfig{file=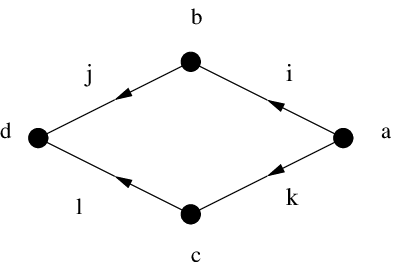,height=3cm}
\caption{\label{f:4cyclepic}}
\end{figure}

 Since there are \(nr\) arrows in the quiver, representations define
 points \((b_i^\rho)\in\mathbb{A}_{\K}^{nr}\).  We write \(\K[z_i^\rho
 : \rho\in G^*, 1\leq i \leq n]\) for the coordinate ring of
 \(\mathbb{A}_{\K}^{nr}\). Our interest lies not with the entire space
 \(\mathbb{A}_{\K}^{nr}\), but with the points \((b_i^\rho)\) of the
 scheme \(Z\) defined by the ideal \[ I = \langle z_j^{\rho\rho_i}
 z_i^{\rho} - z_i^{\rho\rho_j} z_j^{\rho} : \rho \in G^*, 1 \leq i,j
 \leq n\rangle.  \] Thus, we consider only representations
 \((b_i^\rho)\) satisfying the relations 
\begin{equation} \label{e:quiverrelations}
b_j^{\rho\rho_i} b_i^{\rho} = b_i^{\rho\rho_j} b_j^{\rho} \text{ for
 }\rho\in G^* \text{ and }1\leq i,j\leq n
\end{equation}
illustrated in Figure~\ref{f:4cyclepic}.  
These relations arise naturally when quiver representations are
translated into the equivalent language of \emph{\(G\)-constellations}
(see \cite[\S 2]{CMT2}).

 The algebraic torus $\K^*$ acts on each $R_{\rho}$, so \((\K^*)^r\)
 acts diagonally on the vector space \(\oplus_{\rho\in G^*} R_\rho\).
 Hence \(t =(t_\rho)\in (\K^*)^r\) acts on \(b_i^\rho\in
 \Hom(R_{\rho\rho_i},R_{\rho}) = R_{\rho\rho_i}^*\otimes R_{\rho}\) as
\begin{equation} \label{eqn:rtorusaction} t \cdot
 b_i^\rho = t_{\rho\rho_i}^{-1} t_{\rho} b_i^\rho.  \end{equation}
This is a change of basis, with the new basis vector of $R_{\rho}$ set
to be $t_{\rho}^{-1}$ times the old one.  We now describe this action
in terms of a matrix.  Write $\{ \mathbf{e}^{\rho}_i : \rho \in G^{*},
1 \leq i \leq n\}$ for the standard basis of $\Z^{nr}$.  Order the
basis globally into $r$ blocks, one for each \(\rho\in G^*\) beginning
with the trivial representation \(\rho_0\).  Within each block the
elements are listed \(\mathbf{e}^{\rho}_1,\dots,
\mathbf{e}^{\rho}_n\).  Let $B$ be the $r \times (nr)$ matrix with the
column corresponding to $\mathbf{e}_i^{\rho}$ being $\mathbf{e}_{\rho}
- \mathbf{e}_{\rho \rho_i}$, where $\{{\bf e}_\rho : \rho\in G^*\}$ is
the standard basis of $\Z^r$.  Let \(\N B\subset \Z^r\) be the
semigroup generated by the columns of \(B\).  The matrix $B$ is the
vertex-edge incidence matrix of the McKay quiver.  The columns of $B$
encode the weights of the action defined in (\ref{eqn:rtorusaction}).

 \begin{lemma} \label{lemma:ZB} The subsemigroup \(\N B\subset \Z^r\)
   coincides with the sublattice of \(\Z^r\) given by
   \(\emph{\textbf{1}}^\perp := \{(\theta_\rho)\in \Z^r :
   \sum_\rho \theta_\rho = 0\}\).
\end{lemma} 
\begin{proof} The columns of
   $B$ lie in $\textbf{1}^\perp$, so $\N B\subseteq \textbf{1}^\perp$.
   For the opposite inclusion, suppose there is some $\theta \in
   \mathbf{1}^{\perp} \smallsetminus \mathbb N B$.  We may assume that
   $\theta$ has $\sum_{\rho} |\theta_{\rho}|$ minimal for such a
   vector.  Pick $\rho$ with $\theta_{\rho}<0$ and $\rho'$ with
   $\theta_{\rho'}>0$.  The sum of the columns of $B$ corresponding to
   the arrows in a directed path from $\rho$ to $\rho'$ is the vector
   $\textbf{e}_{\rho'} - \textbf{e}_{\rho}$.  Now $\theta':=\theta
   +\mathbf{e}_{\rho}-\mathbf{e}_{\rho'} \in \mathbf{1}^{\perp}$, and
   $\sum_{\rho} |\theta'_{\rho}|$ is smaller, so $\theta' \in \mathbb N
   B$, and thus also $\theta \in \mathbb N B$.
\end{proof}

By Lemma~\ref{lemma:ZB}, the sublattice \(\Z B\subset \Z^r\) generated
by the columns of \(B\) equals the semigroup \(\N B\). The torus
$(\K^*)^r = \Hom(\Z^r,\K^*)$ acts on $\mathbb{A}^{nr}_{\K}$ by formula
(2.3) and induces an action of the torus $T_B:= \Hom(\Z B, \K^*)$.
Note that $T_B = (\K^*)^r/\K^*$.  The ideal \(I\) defining \(Z\) is
invariant under this action, so \(T_B\) acts on \(Z\).

 The \(n\)-dimensional torus \(T^n\) of diagonal matrices with nonzero
 entries acting on \(\mathbb{A}_{\K}^n\) also acts on
 \(\mathbb{A}^{nr}_{\K}\),  where \(s = (s_i)\in T^n\) acts
 on ($b_i^\rho$) as 
 \begin{equation}
   \label{eqn:ntorusaction}
 s\cdot (b_i^\rho) = (s_ib_i^\rho).
 \end{equation}
 Again, the ideal \(I\) defining \(Z\) is invariant under this action,
 so \(T^n\) acts on \(Z\).

 To define moduli spaces of McKay quiver representations we consider
 equivalence classes of quiver representations \((b_i^\rho)\in Z\)
 modulo the \(T_B\)-action.  To construct the quotients we use
 Geometric Invariant Theory (GIT).  For convenience we recall the
 general construction of the GIT quotient of an affine scheme
 $X\subseteq \mathbb{A}^d_{\K}$ by the linear action of an algebraic
 torus $T \cong (\K^{\ast})^s$ (see Dolgachev~\cite{Dolgachev}).
 
 Choose coordinates on $\mathbb{A}^d_{\K}$ to diagonalize the
$T$-action, and identify the character lattice $T^*=\Hom(T,\K^*)$ with
$\Z^s$.  Then for $v = (v_1,\dots, v_d) \in \mathbb{A}^d_{\K}$ and $t
\in T$ there is a character ${\bf a}_i \in T^*=\Z^s$ for which $t
\cdot v_i = t^{{\bf a}_i}v_i = (\prod_{j=1}^s t_i^{(\mathbf{a}_i)_j})
v_i$.  This gives a $\Z^s$-grading of $\K[z_1,\dots,z_d]$ by setting
$\deg(z_i)= {\textbf a_i}$.  Since $T$ acts on $X$, the defining ideal
of $X$ in $\K[z_1,\dots,z_d]$ is homogeneous, and thus $\K[X]$ is
$\mathbb Z^s$-graded.  Line bundles on $X$ are trivial, so a
linearization is a lift of the $T$-action from $X$ to $X\times
\mathbb{A}_{\K}^1$.  If we write $\K[\lambda]$ for the coordinate ring
of $\mathbb{A}^1_{\K}$, then any such lift is determined by $ t\cdot
\lambda = {t}^{-\textbf{b}}\lambda$ for some character $\textbf{b}\in
T^*$.  Then for $f \in \K[X]$ and $j>0$, the function $f \lambda^j \in
\K[X \times \mathbb{A}_{\K}^1]$ is $T$-invariant if and only if
$\deg(f)=j\mathbf{b}$.  Let $\K[X]_{j\textbf{b}}$ be the
$j\mathbf{b}$-graded piece of $\K[X]$. Then the GIT quotient of $X$ by
the action of $T$ linearized by $\textbf{b}$ is
 \[
 X\git_{\!\textbf{b}}T := \Proj \textstyle{\bigoplus_{j\geq 0}}
   \K[X]_{j\textbf{b}}.
 \]
 This scheme is the categorical quotient of the open subscheme of $X$
 consisting of $\textbf{b}$-semistable points.  Recall that a point
 $x\in X$ is $\textbf{b}$-semistable if there exists $j>0$ and $s\in
 \K[X]_{j\textbf{b}}$ such that $s(x)\neq 0$. A
 $\textbf{b}$-semistable point is $\textbf{b}$-stable if the dimension
 of the stabilizer $T_x$ is finite and if there exists some $s\in
 \K[X]_{j\textbf{b}}$ as above for which the $T$-action on the set
 $\{y\in X : s(y)\neq 0\}$ is closed.
 
 Since $\K[X]_\textbf{0}$ is a subalgebra of the graded ring defining
 \(X\git_{\!\textbf{b}}T\), the \(\Proj\) construction induces a
 projective morphism from \(X\git_{\!\textbf{b}}T\) to the quotient
 \(X\git_{\!\textbf{0}}T = \Spec \K[X]^T\) linearized by the trivial
 character.  Moreover, the line bundle \(\owe(1)\) coming from the
 \(\Proj\) construction is relatively ample with respect to this
 morphism.  As is standard in GIT, the quotient linearized by a
 \emph{fractional character} \(\textbf{b}\in T^*\otimes \Q\) is
 defined to be the GIT quotient linearized by any multiple that gives
 an integral character \(j\textbf{b}\in T^*\) (the quotient carries a
 bundle \(\owe(1)\) for which \(\owe(j)\) is relatively ample).  A
 character \(\textbf{b}\in T^*\otimes \Q\) is \emph{generic} if every
 \(\textbf{b}\)-semistable point of \(X\) is in fact
 \(\textbf{b}\)-stable, in which case the categorical quotient
 \(X\git_{\!\textbf{b}}T\) is a geometric quotient.

 We return to the case of interest, where \(T_B = \Hom(\Z B,\K^*)\)
 acts on the affine scheme \(Z\) as in (\ref{eqn:rtorusaction}).  In
 this case, \(T_B^*\otimes \Q\) is the vector space \(\Q B:=\Z
 B\otimes_{\Z} \Q\) generated by the columns of \(B\).  Lemma 2.4
 shows that $\Q B$ is a codimension-one subspace of the vector space
 $\Q^r=\Z^r\otimes_{\Z} \Q$ of fractional characters of the torus
 $(\K^*)^r$.  We introduce the following notation:

 \begin{definition} \label{defn:Theta} The \emph{GIT parameter space}
   is the $\Q$-vector space
 \[
 \Theta := \Q B = \Z B\otimes\otimes_{\Z} \Q = \big{\{}(\theta_{\rho}) \in \Q^r :
 \textstyle{\sum_{\rho \in G^*} \theta_{\rho} = 0}\big{\}}. 
 \]
 For $\theta \in \Theta$, the scheme $\mathcal{M}_\theta:=
 Z\git_\theta T_B$ is the \emph{moduli space of $\theta$-semistable
   McKay quiver representations} of dimension vector \((1,\dots,1)\)
 satisfying the relations~(\ref{e:quiverrelations}).  For generic
 $\theta$,  $\mathcal{M}_\theta$ is the \emph{fine moduli space of $\theta$-stable
   McKay quiver representations}. 
 \end{definition}

 \begin{remark}
   For a finite subgroup \(G\subset \SL(2,\mathbb C)\),
   Kronheimer~\cite{Kronheimer} proved that the geometric quotient
   \(Z\git_\theta T\) coincides with the minimal resolution of
   \(\mathbb{A}^2_{\mathbb C}/G\) for generic \(\theta\in \Theta\).
   The method introduced by Ishii~\cite{Ishii1} extends this result to
   any finite subgroup of \(\GL(2,\K)\).
 \end{remark}


 \section{A distinguished component of $Z$} \label{sec:cohcomp} 

 This section identifies a distinguished irreducible component $V$ of
 the scheme $Z$ introduced in Section~\ref{sec:moduli} and proves that
 \(V\) is a not-necessarily-normal toric variety.
 
Note first that $Z$ need not be irreducible.

\begin{example} \label{ex:reducibleZ} For the action of type
   \(\frac{1}{7}(1,2)\) from Example~\ref{ex:runningexample}, the
   scheme $Z$ is defined by the ideal
$$ I=\langle
z_2^{\rho_1}z_1^{\rho_0} - z_1^{\rho_2}z_2^{\rho_0},
z_2^{\rho_2}z_1^{\rho_1} - z_1^{\rho_3}z_2^{\rho_1},
z_2^{\rho_3}z_1^{\rho_2} - z_1^{\rho_4}z_2^{\rho_2},
z_2^{\rho_4}z_1^{\rho_3} - z_1^{\rho_5}z_2^{\rho_3},$$
$$z_2^{\rho_5}z_1^{\rho_4} - z_1^{\rho_6}z_2^{\rho_4},
z_2^{\rho_6}z_1^{\rho_5} - z_1^{\rho_0}z_2^{\rho_5},
z_2^{\rho_0}z_1^{\rho_6} - z_1^{\rho_1}z_2^{\rho_6}
\rangle.
$$
The ideal $I$ has eight associated primes and hence $Z$ is reducible. 
One of the associated primes is 
$$ J_1 = \langle z_1^{\rho_0}z_2^{\rho_5}-z_1^{\rho_5}z_2^{\rho_6},
z_1^{\rho_6}z_2^{\rho_4}-z_1^{\rho_4}z_2^{\rho_5},z_1^{\rho_5}z_2^{\rho_3}-
z_1^{\rho_3}z_2^{\rho_4},z_1^{\rho_4}z_2^{\rho_2}-z_1^{\rho_2}z_2^{\rho_3},$$
$$z_1^{\rho_3}z_2^{\rho_1}-z_1^{\rho_1}z_2^{\rho_2},z_1^{\rho_6}z_2^{\rho_0}-
z_1^{\rho_1}z_2^{\rho_6},z_1^{\rho_2}z_2^{\rho_0}-z_1^{\rho_0}z_2^{\rho_1},
z_1^{\rho_0}z_1^{\rho_6}z_2^{\rho_3}-z_1^{\rho_3}z_1^{\rho_4}z_2^{\rho_6},$$
$$z_1^{\rho_5}z_1^{\rho_6}z_2^{\rho_2}-z_1^{\rho_2}z_1^{\rho_3}z_2^{\rho_5},
z_1^{\rho_0}z_1^{\rho_6}z_2^{\rho_2}-z_1^{\rho_2}z_1^{\rho_3}z_2^{\rho_6},
z_1^{\rho_5}z_1^{\rho_6}z_2^{\rho_1}-z_1^{\rho_1}z_1^{\rho_2}z_2^{\rho_5},$$
$$z_1^{\rho_4}z_1^{\rho_5}z_2^{\rho_1}-z_1^{\rho_1}z_1^{\rho_2}z_2^{\rho_4},
z_1^{\rho_4}z_1^{\rho_5}z_2^{\rho_0}-z_1^{\rho_0}z_1^{\rho_1}z_2^{\rho_4},
z_1^{\rho_3}z_1^{\rho_4}z_2^{\rho_0}-z_1^{\rho_0}z_1^{\rho_1}z_2^{\rho_3}
\rangle.$$
The distinguished component $V$ of $Z$ is the variety defined by $J_1$.
See Theorem~\ref{t:cohcompeqtns} for details.  Another associated
prime is
$$J_2=\langle z_1^{\rho_1}, z_2^{\rho_1}, z_2^{\rho_2}, z_1^{\rho_5}, z_2^{\rho_5}, z_2^{\rho_6}, z_1^{\rho_4} z_2^{\rho_3} -
z_2^{\rho_4} z_1^{\rho_2}\rangle.$$ The remaining six are obtained by adding \(j\mod
7\) to every raised index in $J_2$, for $j = 1,\dots, 6$.  The
associated primes of $I$ were computed using the computer algebra
package Macaulay~2~\cite{M2}.
\end{example}

  
  We now present a sequence of combinatorial lemmas that enable us to
  define $V$ explicitly in Theorem~\ref{t:cohcompeqtns}. Let $\{
  \mathbf{e}_{\rho}: \rho \in G^{*} \} \cup \{\mathbf{e}_i : 1 \leq i
  \leq n \}$ be the standard basis of $\mathbb Z^{r+n}$. Denote by
  $\pi_n\colon \Z^{r+n}\rightarrow \Z^n$ the projection to the last
  $n$ coordinates.

 \begin{definition} Let $C$ be the $(r+n) \times nr$ matrix whose top
  $r$ rows form the vertex-edge incidence matrix $B$ 
  and whose bottom $n$ rows record the label of the corresponding
  arrow.  Specifically, the column of $C$ indexed by the arrow
  $a_i^{\rho}$ is $C_i^{\rho}:= \mathbf{e}_{\rho}-\mathbf{e}_{\rho
  \rho_i}+\mathbf{e}_i$.  \end{definition}

\begin{example}
   For the action of type \(\frac{1}{7}(1,2)\) from
   Example~\ref{ex:runningexample}, $C$ is the $9 \times 14$ matrix
   shown below and $B$ is the top $7 \times 14$ submatrix.  The third
   column corresponds to the arrow \(a_1^{\rho_1}\) labeled \(i = 1\)
   with tail at \(\rho_2\) and head at \(\rho_1\).  {\Small \[ C=
     \left(\begin{array}{cccccccccccccc}
         1 &  1 &  0 &  0 &  0 &  0 &  0 &  0 &  0 &  0 &  0 & -1 & -1 &
0 \\
         -1 &  0 &  1 &  1 &  0 &  0 &  0 &  0 &  0 &  0 &  0 &  0 &  0 &
-1 \\
         0 & -1 & -1 &  0 &  1 &  1 &  0 &  0 &  0 &  0 &  0 &  0 &  0 &
0 \\
         0 &  0 &  0 & -1 & -1 &  0 &  1 &  1 &  0 &  0 &  0 &  0 &  0 &
0 \\
         0 &  0 &  0 &  0 &  0 & -1 & -1 &  0 &  1 &  1 &  0 &  0 &  0 &
0 \\
         0 &  0 &  0 &  0 &  0 &  0 &  0 & -1 & -1 &  0 &  1 &  1 &  0 &
0 \\
         0 &  0 &  0 &  0 &  0 &  0 &  0 &  0 &  0 & -1 & -1 &  0 &  1 &
1 \\
         1 &  0 &  1 &  0 &  1 &  0 &  1 &  0 &  1 &  0 &  1 &  0 &  1 &
0 \\
         0 & 1 & 0 & 1 & 0 & 1 & 0 & 1 & 0 & 1 & 0 & 1 & 0 & 1
 \end{array}\right)
 \]}
\end{example}

An undirected path (cycle) in the McKay quiver is a path (cycle) in
the underlying undirected graph. Since all arrows in the McKay quiver
are directed, there may be some arrows in the path that are traversed
according to their orientation, and some against.  We write
$-a_i^\rho$ for the arrow $a_i^\rho$ traversed against its
orientation.

\begin{definition} Let $\gamma$ be an undirected path in the McKay
  quiver.  The {\em vector} of $\gamma$, denoted
  $\mathbf{v}(\gamma)\in \mathbb Z^{nr}$, is defined by setting
  $\mathbf{v}(\gamma)_i^{\rho}$ to be the number of times the arrow
  $a_i^{\rho}$ is traversed according to its orientation in the McKay
  quiver minus the number of times it is traversed against its
  orientation.  The {\em type} of a path $\gamma$ is
  $\pi_n(C\mathbf{v}(\gamma)) \in \mathbb Z^n$. The type records the
  number of arrows of each label, where an arrow is counted as
  negative if it is traversed against its orientation. Of particular
  importance are paths of type ${\bf 0}\in \Z^n$.
\end{definition}

\begin{example}
   Consider the McKay quiver for \(\frac{1}{7}(1,2)\) shown in
   Example~\ref{ex:runningexample}. The path $\gamma$ consisting of
   the arrows $a^{\rho_0}_1, a^{\rho_6}_1, -a^{\rho_6}_2,
   -a_1^{\rho_1}$, in that order, has vector $\mathbf{v}(\gamma) =
   (1,0,-1,0,0,0,0,0,0,0,0,0,1,-1)$ and type $(1,-1)$.
 \end{example}

\begin{lemma} \label{l:cyclevector}
  A vector ${\bf u}\in \Z^{nr}$ lies in $\ker_{\mathbb Z}(B):=\{ {\bf u} \in
  \mathbb Z^{nr} \, : \, B{\bf u} = {\bf 0} \}$ if and only if there
  is an undirected cycle $\gamma$ in the McKay quiver with vector
  $\mathbf{v}(\gamma)={\bf u}$.
\end{lemma}

\begin{proof} 
Exercise 38 of Bollob{\'a}s~\cite[II.3]{Bollobas} implies that a vector
$\mathbf{u}$ lies in $\ker_{\mathbb Z}(B)$ if and only if there is a
collection $\{ \gamma_i \}$ of undirected cycles in the McKay quiver
with $\mathbf{u}=\sum_i \mathbf{v}(\gamma_i)$.  Since the McKay quiver
is connected, each $\gamma_i$ can be connected to a base vertex via a
path that is traversed once in the forward direction, and once in the
reverse direction.  The vector of this augmented cycle equals
$\mathbf{v}(\gamma_i)$. Attaching all cycles to this base vertex
produces a single cycle $\gamma$ with vector $\mathbf{v}(\gamma)={\bf
u}$.
\end{proof}

\begin{definition} \label{def:cij}
 For $1 \leq i,j \leq n$ with $i \neq j$, and $\rho \in G^*$, define
 $$\mathbf{c}_{i,j}^{\rho} := \mathbf{e}_i^{\rho} + \mathbf{e}_j^{\rho
   \rho_i} - \mathbf{e}_j^{\rho} - \mathbf{e}_i^{\rho \rho_j} \in
 \mathbb Z^{nr}.$$
\end{definition}

The vectors $\mathbf{c}_{i,j}^{\rho}$ lie in $\ker_{\mathbb Z}(B)$.

\begin{lemma} \label{l:4cyclelem}
  Let $\gamma$ be a path with two adjacent arrows labeled $i$ and $j$,
  with $i \neq j$.  Then there is a path $\gamma'$ satisfying
  $\mathbf{v}(\gamma)-\mathbf{v}(\gamma') = \pm
  \mathbf{c}_{i,j}^{\rho}$ that differs from $\gamma$ only by
  replacing the pair of arrows labeled $i$ and $j$ with a pair labeled
  $j$ and $i$ respectively.
\end{lemma}
\begin{proof}
  The choice of the replacement arrows depends on the orientation of
  the arrows with labels $i$ and $j$, and divides into four cases.

If the arrows labeled $i$ and $j$ meeting at \(\rho\) are both
traversed according to their orientation, replace the arrows
$a_i^{\rho},a_j^{\rho \rho_j^{-1}}$ by $a_j^{\rho \rho_i \rho_j^{-1}},
a_i^{\rho \rho_j^{-1}}$ using $\mathbf{c}_{i,j}^{\rho\rho_j^{-1}}$ as
shown in Figure~\ref{f:arrowlemma}(i).  In this figure the paths go
from right to left, and we replace the path at the top of the diamond
by the one at the bottom.

\begin{figure}[!ht]
\input{switcharrows.pstex_t}
\caption{\label{f:arrowlemma}}
\end{figure}

If the arrows labeled $i$ and $j$ meeting at \(\rho\) are both
traversed against their orientation, replace the arrows $a_i^{\rho
  \rho_i^{-1}},a_j^{\rho}$ by $a_j^{\rho \rho_i^{-1}}, a_i^{\rho
  \rho_j \rho_i^{-1}}$ using $\mathbf{c}_{i,j}^{\rho\rho_i^{-1}}$ as
in Figure~\ref{f:arrowlemma}(ii).
 
If the arrow labeled $i$ is traversed according to its orientation,
but the arrow labeled $j$ against, replace $a_i^{\rho},a_j^{\rho}$ by
$a_j^{\rho \rho_i}, a_i^{\rho \rho_j}$ using $\mathbf{c}_{i,j}^{\rho}$
as in Figure~\ref{f:arrowlemma}(iii).
 
If the arrow labeled $j$ is traversed according to its orientation,
but the arrow labeled $i$ against, replace $a_i^{\rho \rho_i^{-1}},
a_j^{\rho \rho_j^{-1}}$ by $a_j^{\rho \rho_i^{-1} \rho_j^{-1}},
a_i^{\rho \rho_i^{-1} \rho_j^{-1}}$ using $\mathbf{c}_{i,j}^{\rho
  \rho_i^{-1} \rho_j^{-1}}$ as in Figure~\ref{f:arrowlemma}(iv).

 In each case, the orientation of each type of arrow is preserved
 while the labeling of the arrows is switched.
 \end{proof}

The following result, whose proof requires Lemmas~\ref{l:cyclevector}
 and~\ref{l:4cyclelem}, is the key ingredient in
 Theorem~\ref{t:cohcompeqtns}.

\begin{lemma}  \label{c:gensker}
  The set $\{ \mathbf{c}_{i,j}^{\rho} : 1 \leq i,j \leq n, i \neq j,
  \rho \in G^* \}$ of $r \binom{n}{2}$ vectors in $\mathbb Z^{nr}$
  generates the lattice $\ker_{\mathbb Z}(C)=\{ {\bf u} \in \mathbb
  Z^{nr} \, : \, C{\bf u} = {\bf 0} \}$.
\end{lemma}

\begin{proof}
  Let ${\bf u} \in \ker_{\mathbb Z}(C)$.  We need to show that ${\bf
    u}$ can be written as a $\mathbb Z$-linear combination of the
    $\mathbf{c}_{i,j}^{\rho}$.  By Lemma~\ref{l:cyclevector} there is
    a cycle $\gamma$ in the McKay quiver with $\mathbf{v}(\gamma)={\bf
    u}$.  Since ${\bf u} \in \ker_{\mathbb Z}(C)$, $\gamma$ has type
    ${\bf 0}$.  We may assume that any cycle $\gamma'$ consisting of
    fewer arrows than $\gamma$ has $\mathbf{v}(\gamma')$ in the
    integer span of $\mathbf{c}_{i,j}^{\rho}$.
    
    Consider an arrow in $\gamma$ labeled $i$.  By
    Lemma~\ref{l:4cyclelem} we can find a cycle $\gamma'$ with
    $C\mathbf{v}(\gamma')=C\mathbf{v}(\gamma)$ where
    $\mathbf{v}(\gamma)-\mathbf{v}( \gamma')$ is (up to sign) one of
    the $\mathbf{c}_{i,j}^{\rho}$, and the arrow labeled $i$ has been
    switched with an adjacent arrow labeled $j$.  By repeatedly using
    Lemma~\ref{l:4cyclelem} we can move this arrow around until it is
    adjacent to another arrow labeled $i$.  If this new arrow labeled
    $i$ is traversed in the same direction, repeat with this new
    arrow.  Eventually an arrow labeled $i$ traversed in the opposite
    direction must occur, since $\gamma$ has type ${\bf 0}$. But then
    we have a cycle $\gamma''$ with two adjacent arrows with the same
    labels but traversed in opposite directions. Either these two
    arrows labeled $i$ have their heads in a common vertex $\rho$ or
    their tails in a common vertex $\rho$. Since every vertex in the
    McKay quiver has only one incoming and one outgoing edge of the
    same label, the above situation can happen if and only if the
    cycle reaches a vertex using the arrow labeled $i$ and then
    traverses the same arrow in the opposite direction.  Removing this
    cycle of length two produces a cycle $\gamma'''$ satisfying
    $\mathbf{v}(\gamma''')=\mathbf{v}(\gamma'')$. Since $\gamma'''$
    consists of fewer arrows than $\gamma$, the vector ${\bf
      v}(\gamma''')$ lies in the integer span of the
    $\mathbf{c}_{i,j}^{\rho}$ by assumption.  The difference
    $\mathbf{v}(\gamma) - \mathbf{v}(\gamma''')$ also lies in the
    integer span of the $\mathbf{c}_{i,j}^{\rho}$, so
    $\mathbf{v}(\gamma)$ lies in the integer span of the
    $\mathbf{c}_{i,j}^{\rho}$.
\end{proof}

 For ${\textbf u} = (u_i^\rho)\in \mathbb N^{nr}$, we write ${z}^{\bf
 u}$ for the monomial in the polynomial ring $\K[z_i^\rho : \rho\in
 G^*, 1\leq i \leq n]$ that is the product over all $i$ and $\rho$ of
 $z_i^\rho$ raised to the power $u_i^\rho$.  The {\em toric ideal} of
 the matrix $C$ is the ideal
$$I_C := \langle { z}^{\bf u} - {z}^{\bf v} \, : \, {\bf u},
{\bf v} \in \mathbb N^{nr}, \,\, {\bf u}-{\bf v} \in \ker_{\mathbb
Z}(C) \rangle
 $$ in $\K[z_i^\rho : \rho\in G^*, 1\leq i \leq n]$.  As
 Sturmfels~\cite[\S4,13]{GBCP} describes, the toric ideal $I_C$ is a
 prime ideal and defines the not-necessarily-normal affine toric
 variety $\Spec \K[\N C]$.
 
 \begin{theorem} \label{t:cohcompeqtns} There is a unique irreducible
   component of $Z$ that does not lie on any coordinate hyperplane in
   $\mathbb{A}_{\K}^{nr}$.  This component is the affine variety
   $V:=\Spec \K[\N C]$ defined by the toric ideal $I_C$, and is thus
   reduced.
\end{theorem}

\begin{proof}
  By Proposition~\ref{c:gensker} the vectors ${\bf c}_{i,j}^{\rho}$
  generate the lattice $\ker_{\mathbb Z}(C)$. They are precisely the
  differences of exponents of the generators of the binomial ideal $I$
  defining $Z$.  Ho\c{s}ten--Sturmfels~\cite{HostenSturmfels} show
  that if $\{ {\textbf u}_i - {\textbf v}_i : i = 1, \ldots, t \}$
  generates the lattice $\ker_{\Z}(C)$ with $\mathbf{u}_i,
  \mathbf{v}_i \in \mathbb N^{nr}$, and we set $J_C$ to be the ideal
  $\langle { z}^{{\textbf u}_i} - {z}^{{\textbf v}_i} : i = 1, \ldots,
  t \rangle$, then the toric ideal $I_C$ is equal to $(J_C \, : \, (\prod
  z^{\rho}_i)^{\infty})$.  Hence in our case
 $$I_C=(I:(\prod z_i^{\rho})^{\infty}).$$

 The saturation of $I$ by the product of all variables is the ideal of
 the intersection of the components of $Z$ not lying in any coordinate
 hyperplane. Since $I_C$ is prime, its variety is therefore the unique
 component of $Z$ that does not lie entirely in any coordinate
 hyperplane.  This is the variety $V=\Spec \K[\N C]$, since $\K[\N
 C]=\K[\mathbb A_{\K}^{nr}]/I_C$.  Since $I_C$ is prime, $V$ is
 reduced and irreducible.
\end{proof}

\begin{remark}\label{remark:SI}
  Sardo Infirri~\cite[Proposition~5.3]{SI2} claimed that $Z$ is the
  (irreducible) toric variety $\Spec \K[\N C]$. In
  Example~\ref{ex:reducibleZ} we saw that $Z$ may be reducible.  In
  addition, he assumed that $\N C = \Q_{\geq 0}C\cap \Z C$ without
  proof, from which normality of $V = \Spec \K[\N C]$ would follow. We
  show in \cite{CMT2} that $V$ is not always normal.
\end{remark}

\begin{remark} The top $r$ rows of $C$ encode the weights of the $T_B$-action on $\mathbb{A}^{nr}_{\K}$ given by formula (\ref{eqn:rtorusaction}),  and the bottom $n$ rows encode the weights of the $T^n$-action on $\mathbb{A}^{nr}_{\K}$ given by formula (\ref{eqn:ntorusaction}).  Thus,
the action of the dense torus $T_C = \Hom(\Z C,\K^*)$ on $V$ is equal
  to the restriction to $V\subset \mathbb{A}^{nr}_{\K}$ of the
  $T_B\times T^n$-action on $\mathbb{A}^{nr}_{\K}$.
  \end{remark}

 \section{GIT construction of the coherent component}
 \label{sec:toricGIT} 
 
 In this section we show that for all $\theta\in \Theta$, the GIT
 quotient $Y_{\theta}:=V \git_{\theta} T_B$ is a
 not-necessarily-normal toric variety that is obtained from $\mathbb
 A^n_{\K}/G$ by variation of GIT quotient.

 We begin with the construction of \(V\git_{\theta}T_B\).  The action
 of the torus \(T_B\) on \(\mathbb{A}^{nr}_{\K}\) gives a $\mathbb Z
 B$-grading of $\K[\mathbb A_{\K}^{nr}]$ by
 \(\deg(z_i^{\rho})=\mathbf{e}_{\rho}-\mathbf{e}_{\rho \rho_i}\).
 Since $I_C$ is homogeneous in this grading, we obtain a grading of
 the coordinate ring $\K[V]$.  Write $\pi \colon \mathbb Z C
 \rightarrow \mathbb Z B$ for the restriction of the projection $\pi_r
 \colon \mathbb Z^{r+n} \rightarrow \mathbb Z^r$ onto the first $r$
 coordinates.  Since $\K[V]=\K[\mathbb N C]$, the $\mathbb Z
 B$-grading on $\K[V]$ is induced by $\pi$.  In particular, for any
 element \(\theta\in \Theta\) in the GIT parameter space (see
 Definition~\ref{defn:Theta}), the $\theta$-graded piece
 \(\K[V]_{\theta}\) is the \(\K\)-vector space spanned by the
 monomials whose exponents lie in \(\N C\cap \pi^{-1}(\theta)\).  Note
 that this set is nonempty for \(\theta\in \Z B\) by
 Lemma~\ref{lemma:ZB}.  As a result, the categorical quotient of \(V\)
 by the action of \(T_B\) linearized by \(\theta\in \Theta\) is
\begin{equation} \label{eqn:VmodT}
V\git_{\!\theta}T_B = \Proj \textstyle{\bigoplus_{j\geq 0}}
\K[V]_{j\theta} \end{equation} 
where, as before, the definition
for a fractional character 
\(\theta\in \Theta\) is taken to mean
\(V\git_{\!j\theta}T_B\) for some \(j\in \N\) satisfying \(j\theta\in \Z
B\).

Let \(M\subset \Z^n\) be the kernel of the group homomorphism
\(\deg\colon \Z^n\rightarrow G^*\) from equation (\ref{eqn:deg}).  Observe that
\(\mathbb{A}_{\K}^n/G = \Spec \K[\mathbb{A}_{\K}^n]^G = \Spec
\K[\N^n\cap M]\).

 \begin{proposition}
 \label{prop:Yzero}
 The categorical quotient \(V\git_{\!\emph{\textbf{0}}}T_B\) is isomorphic to
 \(\mathbb{A}_{\K}^n/G\).
\end{proposition}

 \begin{proof} For \(\theta = \textbf{0}\) we have
   \(V\git_{\!\textbf{0}}T_B = \Spec \K[V]_{\textbf{0}}\).  Since
   \(\mathbb{A}_{\K}^n/G = \Spec \K[\N^n\cap M]\), the proposition
   follows once we identify \(\N C\cap \ker_\Z(\pi)\) with the semigroup
   \(\N^n\cap M\).
   
   The first step is to show that \(\pi_n\colon \Z^{r+n}\rightarrow
   \Z^n\) induces a lattice isomorphism between \(\ker_\Z(\pi)\) and
   \(M\).  This is equivalent to showing that the respective tori are
   the same.  The lattice \(\ker_\Z(\pi)\) is a sublattice of
   \(\ker_\Z(\pi_r)\) and \(\pi_n(\ker_\Z(\pi))\subseteq
   \pi_n(\ker_\Z(\pi_r)) = \Z^n\). To see that
   \(\pi_n(\ker_\Z(\pi))\subseteq M\), consider \(C\textbf u \in
   \ker_\Z(\pi)\) for $\mathbf{u} \in \mathbb Z^{nr}$.  Since
   \(\textbf u \in \ker_\Z(B)\), Lemma~\ref{l:cyclevector} produces a
   cycle \(\gamma\) in the McKay quiver with \({\textbf v}(\gamma) =
   \textbf u\).  The type of a path from $\rho'$ to $\rho$ in the
   McKay quiver is an element of $\mathbb Z^n$ of degree $\rho^{-1}
   \rho'$, so a path is a cycle if and only if its type lies in $M$.
   Thus \(\pi_n(C\textbf v(\gamma))\) lies in the sublattice
   \(M\subset \Z^n\).  This gives \(\pi_n(C\textbf u)\in M\) as
   claimed. For the opposite inclusion, choose
   \(\mathbf{m}=(m_1,\dots,m_n)\in M\) and construct the path $\gamma$
   in the McKay quiver from any vertex $\rho$ beginning with the
   connected sequence of \(|m_1|\) arrows labeled \(1\), oriented
   according to the sign of $m_1$, followed by the connected sequence
   of \(|m_2|\) arrows labeled \(2\), oriented according to the sign
   of $m_2$, and so on.  This path has type $\mathbf{m}$ by
   construction, so $\gamma$ is a cycle, and thus \({\textbf
     v}(\gamma) \in \ker_\Z(B)\).  This gives \(C{\textbf v}(\gamma)
   \in \ker_\Z(\pi)\), so \(\pi_n(C{\textbf v}(\gamma))\in
   \pi_n(\ker_\Z(\pi))\) and hence \(M\subseteq \pi_n(\ker_\Z(\pi))\).
   Since the restriction of $\pi_n$ to $\ker_\Z(\pi)$ is an
   isomorphism, the claim follows.

   It remains to show that \(\pi_n(\N C\cap \ker_\Z(\pi)) = \N^n\cap M\).
   Since \(\N C\) is generated by vectors of the form \({\bf e}_\rho -
   {\bf e}_{\rho\rho_i}+ {\bf e}_i \), the semigroup \(\pi_n(\N C)\)
   lies in the subsemigroup \(\N^n\subset \Z^n\) generated by the
   elements \(\pi_n({\bf e}_\rho - {\bf e}_{\rho\rho_i}+ {\bf e}_i) =
   {\bf e}_i\) for \(1\leq i\leq n\).  Combining this with the above
   gives \(\pi_n(\N C\cap \ker_\Z(\pi)) \subseteq \N^n\cap M\). To
   establish equality, observe that in the proof of the inclusion
   \(M\subseteq \pi_n(\ker_\Z(\pi))\) described above, if each entry of
   \(\mathbf{m}\in M\) is nonnegative then the path \(\gamma\) has
   \({\textbf v}(\gamma)\in \N^{nr}\).  This gives an element
   \(C{\textbf v}(\gamma)\in \N C\cap \ker_\Z(\pi)\) satisfying
   \(\pi_n({\textbf v}(\gamma))\in \N^n\cap M\).
   \end{proof}

   \begin{remark} \label{remark:standardtorus} The second paragraph of
     the proof of Proposition~\ref{prop:Yzero} shows that the torus
     \(T_C\git_{\!\textbf{0}}T_B = \Spec \K[\Z C]^{T_B} = \Spec
     \K[\ker_\Z(\pi)]\) is isomorphic to the standard torus \(T^n/G =
     \Spec \K[M]\) of the toric variety \(\mathbb{A}_{\K}^n/G\).
   \end{remark}

 \begin{theorem} 
 \label{thm:Ytheta} 
 For $\theta \in \Theta$, set
 \(Y_{\theta}:=V\git_{\theta}T_B\). Then:
 \begin{enumerate}
 \item[\one] $Y_\theta$ is a not-necessarily-normal toric
 variety that admits a projective birational morphism
 $\tau_\theta\colon Y_\theta \rightarrow \mathbb{A}^{n}_{\K}/G$
 obtained by variation of GIT quotient.
 \item[\two] For generic $\theta\in \Theta$, the variety $Y_\theta$ is
 the unique irreducible component of the moduli space $\mathcal
 M_{\theta}$ containing the $T_B$-orbit closures of points of \(Z\cap
 (\K^*)^{nr}\).
 \end{enumerate}
 For generic $\theta\in \Theta$, we call $Y_\theta$ the \emph{coherent
   component} of $\mathcal{M}_\theta$.
\end{theorem} 

\begin{proof} 
  The composition of the canonical projective morphism \(V\git_\theta
  T_B\rightarrow V\git_{\textbf 0}T_B\) described in
  Section~\ref{sec:moduli} with the isomorphism from
  Proposition~\ref{prop:Yzero} gives the projective morphism
  $\tau_\theta$ for all $\theta\in \Theta$. To prove that
  \(\tau_\theta\) is birational we first prove that $T_C
  \git_{\!\theta}T_B$ is isomorphic to $T^n/G$.  Since
  $T_C\git_{\!\textbf{0}}T_B$ is isomorphic to $T^n/G$ by
  Remark~\ref{remark:standardtorus}, we need only show that $T_C
  \git_{\!\theta}T_B$ is isomorphic to $T_C\git_{\!\textbf{0}}T_B$ for
  all $\theta\in \Theta$ or, equivalently, that every point of $T_C$
  is $\theta$-semistable.  To see this, note that each monomial ${\bf
  x}^{\bf u}\in \K[V]_{\theta}$ is nonzero on every point of $T_C$
  because the coordinate entries of every point of $T_C$ are nonzero
  under the given embedding $T_C\subset (\mathbb{\K}^*)^{nr}$.  This
  shows that $T_C$ is $\theta$-semistable, so $\tau_\theta$ is
  birational.
  
  To complete the proof of the first statement we show that $Y_\theta
  = V \git_{\theta}T_B$ is a not-necessarily-normal toric variety.
  The universal property of categorical quotients (see
  \cite[\S0.2]{Mumford}) ensures that \(Y_{\theta}\) inherits
  reducedness and irreducibility from $V$.  Since $T_C\subset V$ is
  dense, and since \(T_C \git_{\theta}T_B\) is nonempty by the above,
  the torus \(T_C \git_{\theta}T_B\) is dense in $Y_\theta$.
  Moreover, the action of $T_C$ on $V$ descends to an action of \(T_C
  \git_{\theta}T_B\) on $Y_\theta$ as required.
  
  To prove the second statement, let $\theta\in \Theta$ be generic, so
  \(V \git_{\theta}T_B\) is a geometric quotient.  The inclusion of
  \(V\) in \(Z\) induces an inclusion of \(Y_{\theta}\) into
  \(\mathcal{M}_{\theta}\).  Since $Y_\theta$ is reduced and
  irreducible, it is a component of \(\mathcal{M}_\theta\) unless
  there is some irreducible component $W\subset Z$ such that
  $W\git_\theta T_B$ contains $V\git_\theta T_B$ as a proper closed
  subscheme.  Let $V^s_\theta$ and $W^s_\theta$ denote the loci of
  $\theta$-stable points of $V$ and $W$ respectively, and write
  $p\colon W^s_\theta\rightarrow W\git_\theta T_B$ for the natural
  quotient map.  Since $V\git_\theta T_B\subseteq W\git_\theta T_B$,
  and since the fibers of $p$ are closed $T_B$-orbits, we obtain
  $V^s_\theta\subseteq W^s_\theta$.  This implies that $W^s_\theta$
  contains the torus $T_C$, because every point of $T_C\subset V$ is
  $\theta$-stable by the first paragraph above.  Then the irreducible
  component $W$ of the scheme $Z$ contains points of $Z \cap
  (\K^*)^{nr}$, which contradicts Theorem~\ref{t:cohcompeqtns} since
  $W\neq V$.  Thus, $Y_\theta$ is an irreducible component of
  $\mathcal{M}_\theta$.
 \end{proof}

 \begin{remark} \label{remark:coherent} The adjective \emph{coherent}
   comes from the theory of the toric Hilbert scheme introduced by
   Peeva--Stillman~\cite{PeevaStillman}.
 \end{remark}

  \begin{remark}
    For finite abelian \(G\subset \SL(n,\K)\) and $\theta\in \Theta$
    with \(n\leq 3\), the projective birational morphism
    \(\tau_\theta\colon Y_\theta \rightarrow \mathbb{A}_{\K}^n/G\)
    from Theorem~\ref{thm:Ytheta} is a crepant resolution by
    Kronheimer~\cite{Kronheimer} and
    Bridgeland--King--Reid~\cite{BKR}.  Their results hold without the
    abelian assumption.
\end{remark}

Theorem~\ref{thm:Ytheta} suggests the following conjecture that does
not require the abelian assumption on $G$.  For any finite subgroup
$G\subset \GL(n,\K)$, moduli spaces of $\theta$-stable quiver
representations are constructed as GIT quotients $\mathcal{M}_\theta =
Z\git_\theta H$, where $Z$ is an affine scheme, $H$ is an algebraic
group and \(\theta\in H^*\otimes \Q\) is a fractional character (see
Craw--Ishii~\cite[\S2]{CrawIshii}).  The $G$-orbit of the point
$(1,\dots,1)\in \mathbb{A}^n_{\K}$ defines a quiver representation
$z\in Z$ that gives a point $[z]\in \mathcal{M}_\theta$ for all
$\theta\in \Theta$.  The algebraic torus $(\K^*)^n$ also acts on $Z$.

\begin{conjecture}
  Let $V$ denote the subscheme of $Z$ obtained as the closure of the
  $(\K^*)^n\times H$-orbit of $z\in Z$. For generic $\theta \in
  H^*\otimes \Q$, the GIT quotient \(Y_{\theta}:=V\git_{\theta}H\) is
  a reduced irreducible component of $\mathcal M_{\theta}$ that admits
  a projective birational morphism \(Y_\theta \rightarrow Y_{\bf 0}
  \cong \mathbb{A}_{\K}^n/G\) obtained by variation of GIT quotient.
\end{conjecture}

 \section{The \protect$G\protect$-Hilbert scheme}
 \label{sec:hilbg}
 Applying Theorem~\ref{thm:Ytheta} to the special case where $\mathcal
 M_\theta \cong \ghilb$ provides a simple construction of Nakamura's
 $G$-Hilbert scheme.  Assume that $G\subset \GL(n,\K)$ is abelian.
 Let $S$ denote the coordinate ring $\K[x_1,\dots, x_n]$ of
 $\mathbb{A}^n_{\K}$.
 
 We first recall the $G$-Hilbert scheme. The literature contains
 two inequivalent definitions as follows.  The first, denoted
 \(\ghilb\), is the fine moduli space of ideals \(J\subseteq S \)
 defining \(G\)-invariant subschemes \(Z(J)\subseteq
 \mathbb{A}^{n}_{\K}\) whose coordinate rings \(S/J\) are isomorphic to
 $\K G$ as a \(\K G\)-module. The ideal \(J\), or the scheme \(Z(J)\),
 is called a \emph{\(G\)-cluster}, and we write \([J]\in \ghilb\).

 To define the second (and original) version, denoted \(\hilbg\),
 observe that the \(G\)-orbit of any point \(p\in T^n\subseteq
 \mathbb{A}^n_{\K}\) consists of \(r\) distinct points permuted
 transitively by $G$.  These orbits define points $[G\cdot p]$ in the
 \(G\)-fixed locus $(\hilb^r(\mathbb{A}^n_{\K}))^G$ in the Hilbert
 scheme of \(r\) points in \(\mathbb{A}^n_{\K}\).  Every such point
 lies in a unique irreducible component of
 $(\hilb^r(\mathbb{A}^n_{\K}))^G$ that we denote \(\hilbg\).  To see
 this, note that an infinitesimal \(G\)-equivariant deformation of
 \([G\cdot p]\in \hilb^r(\mathbb{A}^n_{\K})\) with \(p\in T^n\)
 deforms the \(r\) distinct points in \(T^n\) that support the orbit.
 The resulting subscheme \(Z'\subseteq \mathbb{A}^n_{\K}\) is
 supported on \(r\) distinct points in \(T^n\) and, since the
 deformation was \(G\)-equivariant, we obtain \(Z' = G\cdot p'\) for
 some point \(p'\in T^n\).  This shows that \(G\)-clusters of the form
 \([G\cdot p]\) for \(p\in T^n\) are open in some union of components
 of $(\hilb^r(\mathbb{A}^n_{\K}))^G$.  There is only one component
 since for any two such \(G\)-clusters \([G\cdot p], [G\cdot p']\in
 \ghilb\), one can construct a morphism \(\mathbb{A}^1\rightarrow
 \ghilb\) whose image contains both \([G\cdot p]\) and \([G\cdot
 p']\). This shows that $\hilbg$ is well-defined.

   Ito--Nakajima~\cite[\S2]{ItoNakajima} proved that $\ghilb$ is a
   union of connected components of $(\hilb^r(\mathbb{A}^n_{\K}))^G$.
   This also follows from
   Haiman--Sturmfels~\cite[Proposition~1.5]{HaimanSturmfels}.  The
   \(G\)-orbits \(G\cdot p\subseteq \mathbb{A}^{n}_{\K}\) defined by
   points \(p\in T^n\) are $G$-clusters, so the component of
   $(\hilb^r(\mathbb{A}^n_{\K}))^G$ containing points of the form
   \([G\cdot p]\) for \(p\in T^n\) is a component of \(\ghilb\).  This
   is \(\hilbg\) by definition.  In particular, $\hilbg \subseteq
   \ghilb$.

 \begin{remark}
   The original definition of \(\hilbg\) is due to
   Ito--Nakamura~\cite{ItoNakamura} and further studied by
   Nakamura~\cite{Nakamura}. The moduli definition \(\ghilb\) is
   due to Reid~\cite{Reid2} and is the version of the
   \(G\)-Hilbert scheme adopted by \cite{ItoNakajima, BKR, Ishii1,
     HaimanSturmfels}.  For a finite subgroup \(G\) of \(\GL(2,\K)\)
   or \(\SL(3,\K)\), it is known that \(\ghilb\) is smooth and
   connected (see \cite{Ishii1, BKR}), hence \(\ghilb \cong \hilbg\).
 \end{remark}

 These distinct versions of the \(G\)-Hilbert scheme can be
 constructed simultaneously via moduli of quiver representations as
 follows.  Also, Nakamura~\cite{Nakamura} asserted that \(\hilbg\) is
 endowed with the reduced scheme structure, but here we show that this
 follows naturally from the definitions and Theorem~\ref{thm:Ytheta}.

 \begin{proposition}
   \label{prop:hilbg}
   Let \(\theta\in \Theta\) satisfy $\theta_{\rho_0} < 0$ and
   \(\theta_\rho > 0\) for \(\rho\neq \rho_0\).  Then there is an
   isomorphism \(\mathcal{M}_\theta \cong \ghilb\) that induces an
   isomorphism \(Y_\theta \cong\hilbg\) by restriction.  In particular, the
   scheme \(\hilbg\) is reduced.
 \end{proposition}
 \begin{proof}
   Ito--Nakajima~\cite[\S3]{ItoNakajima} observed that there is a
   unique chamber \(C_+\subseteq \Theta\) in the GIT parameter space
   containing parameters \(\{\theta \in \Theta \st
   \theta_\rho > 0 \text{ for } \rho\neq \rho_0\}\) such that
   \(\mathcal{M}_\theta \cong \ghilb\) for all \(\theta\in C_+\).  To
   complete the proof of the first statement it remains to show that
   this isomorphism identifies the coherent component
   \(Y_{\theta}\subseteq \mathcal{M}_\theta\) with
   \(\hilbg\subseteq\ghilb\).  This follows from the proof of
   Theorem~\ref{thm:Ytheta}, where it is shown that the standard torus
   of \(Y_{\theta}\) is isomorphic to the standard torus
   \(T^n/G\subseteq \mathbb{A}^{n}_{\K}/G\) parametrizing \(G\)-orbits
   \(G\cdot p\) for \(p\in T^n\).  The final statement follows since
   \(Y_\theta\) is reduced by Theorem~\ref{thm:Ytheta}.
 \end{proof}

 \begin{remark} 
   The notation $\mathcal M_\theta \cong \ghilb$ means that not only
   are the underlying schemes isomorphic but, in addition, the
   tautological bundles on both $\mathcal M_\theta$ and $\ghilb$
   induced by the moduli constructions are also isomorphic as
   $G$-equivariant locally free sheaves.
\end{remark}

 \begin{remark}
   Proposition~\ref{prop:hilbg} provides a direct GIT construction of
   the irreducible scheme \(\hilbg \cong V\git_\theta T_B\) that
   avoids the Hilbert scheme of \(r\)-points in \(\mathbb{A}^n_{\K}\),
   and shows that \(\hilbg\) may be obtained from \(\mathbb{A}^n/G\)
   by variation of GIT quotient.  The Hilbert scheme of \(r\)-points
   in \(\mathbb{A}^n_{\K}\) is in general much more singular than
   anything needed for $\hilbg$.
 \end{remark}

\section{Effective computation of the fan of \protect\(Y_\theta\protect\)}
\label{sec:toricfan}
The GIT construction of $Y_\theta$ allows an explicit description of
the fan of the not-necessarily-normal toric variety $Y_{\theta}$.

\begin{definition} \label{d:ptheta}
Let \(P_\theta\) denote the convex hull of the set \(\pi_n(\N C\cap
\pi^{-1}(\theta))\) in the vector space
\(\pi_n(\pi^{-1}(\theta))\otimes_{\Z}\Q\).  Since the lattices
\(\pi^{-1}(\theta)\) and \(\ker_\Z(\pi)\) are isomorphic, the proof of
Proposition~\ref{prop:Yzero} gives \(\pi_n(\pi^{-1}(\theta))\cong M\).
As a result, we regard \(P_\theta\) as a polyhedron in
\(M\otimes_{\Z}\Q\) for all \(\theta\in \Theta\).

Let $F$ be a face of $P_{\theta}$. The \emph{inner normal cone}
$\mathcal N_{P_{\theta}}(F)$ of $P_{\theta}$ at $F$ is the set of all
${\textbf y} \in M^{\vee} \otimes \mathbb Q$ such that the linear
functional ${\textbf y}$ is minimized over $P_{\theta}$ at $F$. The
{\em inner normal fan} \(\mathcal{N}(P_{\theta})\) of $P_{\theta}$ is
the polyhedral fan whose cones are $\{\mathcal N_{P_{\theta}}(F)\}$ as
$F$ varies over the faces of $P_{\theta}$.
\end{definition}

We now describe $Y_{\theta}$ as a not-necessarily-normal toric variety
in terms of a fan (plus extra data) following
Thompson~\cite{Thompson}.  Replace the parameter \(\theta\) by a
positive multiple if necessary to ensure that the graded ring
\(\oplus_{j\geq 0} \K[V]_{j\theta}\) defining \(Y_\theta\) is
generated in degree one.  Then \(Y_\theta\) is covered by charts of
the form \(\Spec ((\oplus_{j\geq 0} (\K[V]_{j\theta})_{t})_0)\), where
\(t\) varies over the generators of some ideal with radical the
irrelevant ideal \(\oplus_{j> 0} \K[V]_{j\theta}\), and where
$(\oplus_{j\geq 0} (\K[V]_{j\theta})_{t})_0$ denotes the degree-zero
piece of the localization of $\oplus_{j\geq 0} \K[V]_{j\theta}$ at
$t$.  We choose as a generating set the set of vertices of the convex
lattice polyhedron \(P_\theta\).  Let $\mathbf{m}$ be a vertex of
$P_{\theta}$, with $\sigma=\mathcal N_{P_{\theta}}(\mathbf{m})$, and
let \(A_{\sigma}\) be the subsemigroup of \(M\) given by
$A_{\sigma}:=\mathbb N\langle \mathbf{p}-\mathbf{m} : \mathbf{p} \in
P_\theta\cap M\rangle$. Then the affine charts covering $Y_{\theta}$
are of the form \(\Spec \K[A_\sigma]\), where $\sigma$ varies over the
normal cones of the vertices of $P_{\theta}$.

  If \(Y_\theta\) is normal then the semigroups \(A_{\sigma}\) can be
  written as \(A_{\sigma} = \sigma^\vee\cap M\).  It follows that
  \(Y_\theta\) is the toric variety with fan
  \(\mathcal{N}(P_\theta)\subseteq M^{\vee} \otimes \Q\).  Otherwise,
  \(A_{\sigma} \subsetneq \sigma^\vee\cap M\) in general and
  \(Y_\theta\) is described by a fan \((\Delta, \mathcal{S})\) in the
  sense of Thompson~\cite{Thompson} as follows.  Consider the set
  \(\Delta:= \{\sigma\in \mathcal{N}(P_\theta)\}\) as the topological
  space whose open sets are the subfans of \(\mathcal{N}(P_\theta)\),
  where each cone \(\sigma\in \mathcal{N}(P_\theta)\) is regarded as
  the fan consisting of the faces of \(\sigma\).  We define a sheaf
  $\mathcal S$ of semigroups on \(\Delta\) by first setting
  \(\Gamma(\sigma, \mathcal{S}\vert_{\sigma}) := A_{\sigma}\) for
  \(\sigma = \mathcal{N}_{P_\theta}(m)\).  More generally, if $\tau$
  is a face of $\sigma=\mathcal{N}_{P_{\theta}}(m)$, then
  \(\Gamma(\tau, \mathcal{S}\vert_{\tau}) := A_{\tau}\) where
  $A_{\tau}=A_{\sigma} + \mathbb Z \mathbf{u}$ for $\tau = \sigma
  \cap \mathbf{u}^{\perp}$ (compare \cite[Lemma 1.3]{Fulton}).  Then
  the pair \((\Delta, \mathcal{S})\) defines the nonnormal toric
  variety \(Y_\theta\) as in Thompson~\cite{Thompson}.  The point is
  simply that \(Y_\theta\) is built from the local charts \(\Spec
  \K[A_{\sigma}]\) rather than the normal varieties \(\Spec
  \K[\sigma^\vee\cap M]\).  In particular, the normalization of
  $Y_{\theta}$ is obtained by replacing each $A_{\sigma}$ by its
  normalization $\sigma^{\vee} \cap M$.

 \begin{corollary} 
 \label{cor:fan-of-normalization}
 The normalization $\widetilde{Y_\theta}$ of \(Y_\theta\) is the toric
 variety whose fan is $\mathcal N(P_{\theta})$. Furthermore, the fan
 $\mathcal N(P_{\theta})$ is supported on the cone $(\Q^n_{\geq 0})^\vee$.
 \end{corollary}
 \begin{proof}
   It remains to prove the second statement.  The morphism from
   Theorem~\ref{thm:Ytheta} induces a projective, birational toric
   morphism $\widetilde{Y_\theta}\rightarrow \mathbb{A}^n_{\K}/G$, so
   the support of $\mathcal N(P_{\theta})$ equals the support of the
   cone $(\Q^n_{\geq 0})^\vee$ in $M^\vee\otimes \Q$ defining
   $\mathbb{A}^n_{\K}/G$.
\end{proof}

We call $\mathcal N(P_{\theta})$ the \emph{fan of} $Y_{\theta}$.

\begin{remark} Sardo Infirri~\cite[Theorem~5.5]{SI2} claimed that
  \(\mathcal{M}_\theta\) is the toric variety with fan $\mathcal
  N(P_{\theta})$ for all \(\theta\in \Theta\).  We show in
  \cite[Examples~4.12, 5.6]{CMT2} that $\mathcal{M}_\theta$ is not a
  normal toric variety in general.

\end{remark}

 \begin{example}
 \label{ex:weightoneaction}
  Let $G \cong \mathbb Z/r\Z \subset \GL(n,\K)$ be generated by
   $\diag(\omega, \dots, \omega)$, where $\omega$ is a primitive $r$th
   root of unity, so $\mathbb A_{\K}^n/G$ is of type
   $\frac{1}{r}(1,\dots,1)$.  The lattice $M^{\vee}$ is generated by
   the standard basis vectors plus the vector
   $\frac{1}{r}(1,\dots,1)$. As a toric variety, $\mathbb A_{\K}^n/G$
   corresponds to the rational cone generated by the standard basis
   vectors in $M^{\vee}\otimes\Q$.  It is well known that a resolution
   $Y\rightarrow \mathbb{A}^n_{\K}/G$ is obtained by adding the ray
   generated by $\frac{1}{r}(1,\dots,1) \in M^{\vee}$ and taking the
   stellar subdivision (see, for example, \cite[Example 4.7]{Craw}).
   This resolution is crepant if and only if \(n=r\).
   
   We now show that $Y$ is isomorphic as a variety to $Y_{\theta}$ for
   all $\theta \in \Theta \smallsetminus \{ 0 \}$.  For $\theta \in
   \Theta$, the polyhedron $P_{\theta}$ is equal to $\pi_n\{ C {\bf u} \,:\, B{\bf
     u} = \theta, {\bf u} \geq {\bf 0} \}$.  Since $\pi_n (C {\bf u})
   = \sum u_i^\rho {\bf e}_i$, the minimum $1$-norm of a vector in
   $P_{\theta}$ is $d_\theta := \textup{min} \{ \sum u_i^\rho \,:\,
   {\bf u} \geq {\bf 0}, B {\bf u} = \theta \}$, so $P_{\theta}
   \subseteq \{ \mathbf{y} \in \mathbb Q^n_{\geq 0} : \sum_i y_i \geq
   d_{\theta} \}$.  Since $B$ is a unimodular matrix, we know that
   $d_{\theta} \in \mathbb N$.  We claim that this inclusion is
   equality.  First, we show that the vertices of the right hand side
   lie in $P_\theta$.  Indeed, by the unimodularity of $B$ there
   exists ${\bf y} \in P_\theta \cap M$ whose 1-norm is $d_{\theta}$.
   By replacing each ${\bf e}_i^\rho$ by ${\bf e}_1^\rho$, any vector
   ${\bf u} = \sum u_i^\rho {\bf e}_i^\rho \geq {\bf 0}$ satisfying $C
   {\bf u} = (\theta,{\bf y})$ determines a vector ${\bf u}' := \sum
   u_i^\rho {\bf e}_1^\rho$ of type $d_\theta {\bf e}_1^\rho$.  Since
   the McKay quiver of $G$ is a cycle with $n$ arrows connecting
   adjacent vertices, we have $B{\bf u}' = B{\bf u}$ and hence $C {\bf
     u}' = (\theta, d_\theta {\bf e}_1)$. This shows that the vertex
   $d_\theta {\bf e}_1$ and, similarly, any vertex $d_\theta {\bf
     e}_i$, lies in $P_\theta$. To show that any point on a facet of
   $\{ \mathbf{y} \in \mathbb Q^n_{\geq 0} : \sum_i y_i \geq
   d_{\theta} \}$ other than $\{ \mathbf{y} \in \mathbb Q^n_{\geq 0} :
   \sum_i y_i = d_{\theta} \}$ also lies in $P_\theta$, note that if
   $\gamma$ is a cycle in the McKay quiver of type $r{\bf e}_i$ then
   $C({\bf u}' + j{\bf v}(\gamma)) = (\theta,(d_\theta + jr){\bf
     e}_i)$ for all $j > 0$.  This proves that $P_{\theta} = \{
   \mathbf{y} \in \mathbb Q^n_{\geq 0} : \sum_i y_i \geq d_{\theta}
   \}$. Note that $d_{\theta}=0$ if and only if $\theta=0$.  When
   $d_{\theta}>0$, the normal fan of $P_\theta$ is the fan of the
   resolution $Y\rightarrow \mathbb{A}^n_{\K}/G$ defined above.  In
   addition, it is easy to see that the semigroups satisfy
   $A_{\sigma}=\sigma^{\vee} \cap M$ for all top-dimensional cones
   $\sigma$ in $\mathcal N(P_{\theta})$.  Thus $Y_{\theta}$ is normal,
   so $Y\cong Y_\theta$ for any $\theta \in \Theta \smallsetminus \{ 0
   \}$.  It follows from Theorem~\ref{thm:Ytheta} that for $\theta\in
   \Theta$ generic, the toric resolution $Y$ of $\mathbb{A}^n_{\K}/G$
   is isomorphic to the coherent component of the moduli space
   $\mathcal{M}_\theta$.
\end{example}

To emphasize that the previous results are explicit, we now give an
algorithm to compute the fan of $Y_\theta$. Every polyhedron $P$ can
be written as the Minkowski sum of a polytope $Q$ and a polyhedral
cone $K$. By the {\em generator representation} of $P$ we mean the
pair of lists $L_1, L_2$, where $L_1$ consists of the vertices of $Q$,
and $L_2$ consists of the generators of $K$.  The software package
PORTA \cite{porta} converts between the inequality and generator
representation of a polyhedron.  Recall that the top $r \times nr$
submatrix of the matrix $C$ is the vertex-edge incidence matrix $B$ of
the McKay quiver. Let $D$ be the bottom $n \times nr$ submatrix of
$C$.

 \begin{algorithm} To compute the fan of $Y_{\theta}$.
 \end{algorithm}

 \noindent{\bf Input:} The GIT parameter $\theta \in \Theta$ and  
 the matrix $C$.

\begin{enumerate} 
\item Compute a generator representation for the polyhedron 
$\{ {\bf u} \in \mathbb Q_{\geq 0}^{nr} \,:\, B{\bf u} = \theta \}$ and obtain
the lists $L_1$ and $L_2$.
\item For $i = 1,2$, replace the list $L_i$ with the list 
  $DL_i := \{ D{\bf u} \,:\, {\bf u} \in L_i \}$.
\item Compute the inequality description of the polyhedron
  $P_{\theta}$ obtained as the sum of the polytope $\textup{conv}
  (DL_1)$ and the cone generated by $DL_2$.
\item The output of PORTA contains a table of incidences between
  inequalities of $P_{\theta}$ and vertices of $P_{\theta}$. The
  normal fan at a vertex is generated by the negatives of normal
  vectors on the left hand side of those inequalities that hold at
  equality at the vertex.
\end{enumerate}

\noindent{\bf Output:} The inner normal fan of $P_{\theta}$ as a
collection of sets of generators of normal cones at the vertices of
$P_{\theta}$.

\begin{proof}[Proof of Correctness]
  The fan of ${Y_\theta}$ is the inner normal fan of the lattice
  polyhedron $P_\theta = \conv(\pi_n(\N C\cap \pi^{-1}(\theta)))$. The
  set of lattice points $\N C\cap \pi^{-1}(\theta) = C (\{{\bf u} \in
  \mathbb N^{nr} \, : \, B {\bf u} = \theta \})$ is the image of
  $\{{\bf u} \in \mathbb N^{nr} \, : \, B {\bf u} = \theta \}$ under
  the linear map $C \colon\mathbb Q^{nr} \rightarrow \mathbb Q^{r+n}$.
  Further, $\pi_n(\N C\cap \pi^{-1}(\theta)) = D (\{{\bf u} \in
  \mathbb N^{nr} \, : \, B {\bf u} = \theta \})$ and hence $P_\theta =
  D (\conv (\{{\bf u} \in \mathbb N^{nr} \, : \, B {\bf u} = \theta
  \}))$ since the operation of taking convex hulls commutes with
  linear maps. Since the matrix $B$ is totally unimodular (see
  Schrijver~\cite[p274, Example 2]{Schrijver}), we have $\conv (\{{\bf
    u} \in \mathbb N^{nr} \, : \, B {\bf u} = \theta \}) = \{ {\bf u}
  \in \mathbb Q^{nr}_{\geq 0} \, : \, B {\bf u} = \theta \}$.  Hence
  $P_\theta = D (\{ {\bf u} \in \mathbb Q^{nr}_{\geq 0} \, : \, B {\bf
    u} = \theta \})$. This justifies steps (1) and (2) of the
  algorithm.  Step (4) extracts the normal fan of $P_{\theta}$ by
  computing the normal cone at each vertex of $P_{\theta}$.
\end{proof}

 \begin{example}
 \label{ex:128} 
   Consider the action of type $\frac{1}{11}(1,2,8)$, so $G\subset
   \SL(3,\K)$ is the cyclic group of order $11$ with generator
   $\diag(\omega, \omega^2, \omega^8)$ where $\omega$ is a primitive
   11th root of unity. We compute the fan $\mathcal{N}(P_{\theta})$ for
   $\theta = (1,1,1,1,-7,-9,1,1,1,8,1).$

   In this example, $r=11$ and $n=3$.  The matrix $C$ is a $14 \times
   33$ matrix, and the polyhedron $\{{\bf u} \in \mathbb Q_{\geq
   0}^{nr} \,:\, B{\bf u} = \theta \}$ is the Minkowski sum of a cone
   generated by $630$ vectors and a polytope with $17581$
   vertices. This computation was done using PORTA~\cite{porta}. After
   multiplying these lists of vectors by $D$, we compute $P_{\theta}
   \subset \mathbb Q^3$ as the sum of the convex hull of $DL_1$ and
   the cone generated by $DL_2$. Then $P_{\theta}$ is the sum of the
   cone generated by $(0,0,1), (0,1,0), (1,0,0)$ and the convex hull
   of $(0,0,78)$, $(0,21,15)$, $(0,26,11)$, $(0,70,0)$, $(22,0,23)$,
   $(96,0,0)$, $(4,0,50)$, $(4,9,23)$, $(4,46,0)$, $(72,0,3)$,
   $(4,34,3)$. Equivalently, $P_{\theta}$ is described by the
   irredundant inequalities listed in Table~\ref{t:portacalc1}.
 \begin{table}[!ht]
\mbox{\subfigure[]{
$\begin{array}{llll}
(1) &  -2x_1-4x_2-5x_3 & \leq & -159\\
(2) & -3x_1-6x_2-2x_3  & \leq & -112\\
(3) &  - x_1-2x_2-8x_3 & \leq &  -96\\
(4) &  -7x_1-3x_2- x_3 & \leq &  -78\\
(5) &  -6x_1- x_2-4x_3 & \leq &  -70\\
(6) &  - x_1         & \leq &    0\\
(7) &      - x_2     & \leq &    0\\
(8) &          - x_3 & \leq &    0\\
\end{array}$\label{t:portacalc1}}
\hspace{1.5cm} 
\subfigure[]{
$\begin{array}{|l|l|}
\hline
(0, 0,78) & \{4,6,7\} \\
(0,21,15) & \{1,4,6\} \\
(0,26,11) & \{1,5,6\} \\
(0,70, 0) & \{5,6,8\} \\
(22,0,23) & \{1,2,7\} \\
(96,0 ,0) & \{3,7,8\} \\
(4, 0,50) & \{2,4,7\} \\
(4, 9,23) & \{1,2,4\} \\
(4,46 ,0) & \{3,5,8\} \\
(72 ,0 ,3) & \{1,3,7\} \\
(4,34 ,3)  & \{1,3,5\} \\
\hline
\end{array}$
\label{t:portacalc2}}}
 \caption{(a) inequalities; (b) output from PORTA}
 \label{t:portacalc}
\end{table}

To obtain the normal cones of $P_{\theta}$ we calculate the
inequalities that hold at equality at which vertex. This information
is carried in the strong validity table from PORTA at the end of the
computation. We condense this information in Table~\ref{t:portacalc2}.
The first line of this table means that the normal cone at the vertex
$(0,0,78)$ of $P_{\theta}$ is generated by the negatives of the
coefficient vectors of the linear forms on the left hand side of
inequalities $(4),(6)$ and $(7)$, in this case $(7,3,1),(1,0,0)$ and
$(0,1,0)$.  A cross-section of the fan $\mathcal{N}(P_\theta)$ is
shown in Figure~\ref{fig:128}.  The rays of $\mathcal{N}(P_\theta)$
are labeled 1 through 8 according to Table~\ref{t:portacalc1} so, for
example, ray 6 is generated by $(1,0,0)$ and ray 3 is generated by
$\frac{1}{11}(1,2,8)$.
 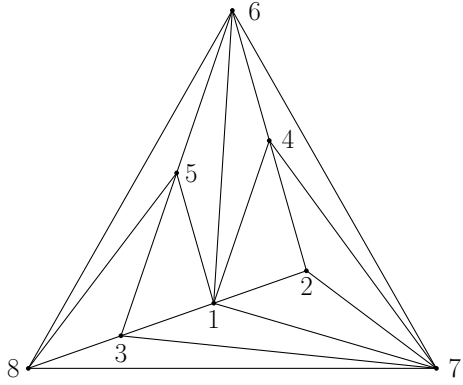
\begin{figure}[!ht]
 \centering
 \psset{unit=1cm}
 \resizebox{5cm}{5cm}{%
 \begin{pspicture}(10,10)
 \psline{-}(0,0)(11,0)
 \psline{-}(11,0)(5.5,9.526)
 \psline{-}(5.5,9.526)(0,0)
 \psdots(0,0)(11,0)(5.5,9.526)
 \rput(6.1,9.5){\LARGE{6}}
 \rput(11.5,0){\LARGE{7}}
 \rput(-0.4,0){\LARGE{8}}
 \psdots(2.5,0.866)(5,1.732)(7.5,2.598)(4,5.196)(6.5,6.062)
 \rput(2.5,0.4){\LARGE{3}}
 \rput(5,1.3){\LARGE{1}}
 \rput(7.5,2.2){\LARGE{2}}
 \rput(4.4,5.2){\LARGE{5}}
 \rput(7,6.1){\LARGE{4}}
 \psline(0,0)(2.5,0.866)
 \psline(0,0)(4,5.196)
 \psline(6.5,6.062)(11,0)
 \psline(7.5,2.598)(11,0)
 \psline(5,1.732)(11,0)
 \psline(2.5,0.866)(11,0)
 \psline(4,5.196)(5.5,9.526)
 \psline(5.5,9.526)(6.5,6.062)
 \psline(2.5,0.866)(5,1.732)
 \psline(5,1.732)(7.5,2.598)
 \psline(2.5,0.866)(4,5.196)
 \psline(6.5,6.062)(7.5,2.598)
 \psline(4,5.196)(5,1.732)
 \psline(5,1.732)(6.5,6.062)
 \psline(5,1.732)(5.5,9.526) 
 \end{pspicture}
 }
 \caption{The fan of $Y_\theta$ for $\theta = (1,1,1,1,-7,-9,1,1,1,8,1)$}
 \label{fig:128}
 \end{figure}

 We note that while we have shown all steps in the algorithm
 explicitly in this example, this procedure can be completely
 automated.
  \end{example}

 \begin{remark}
   The calculation of the toric fan defining $Y_\theta$ as in
   Figure~\ref{fig:128} was originally done by hand by
   Craw~\cite[\S5.8.2]{Crawthesis} using a lengthy and somewhat
   speculative method.  The method implemented here generalizes easily
   to calculate the fan for significantly more involved examples, and
   is automated.
 \end{remark}

\section{Distinguished McKay quiver representations}
 \label{sec:distinguished}

 In this section we calculate the distinguished $\theta$-semistable
 quiver representations that define points on $Y_\theta$.  For generic
 $\theta\in \Theta$, these representations encode the restriction to
 $Y_\theta$ of the universal quiver representation on the fine moduli
 space $\mathcal{M}_\theta$ (see Craw--Ishii~\cite[\S2]{CrawIshii}).

 Recall from Theorem~\ref{t:cohcompeqtns} that $V=\Spec \K[\N C]$,
 where $C$ is the matrix with columns
 $C_i^{\rho}=\mathbf{e}_{\rho}-\mathbf{e}_{\rho \rho_i}+\mathbf{e}_i$.
 The standard torus of $V$ is $T_C$, and, as for normal toric
 varieties, the orbits of the $T_C$-action correspond one-to-one to
 faces of the rational polyhedral cone $P:= \Q_{\geq 0} C \subseteq
 \Q^{r+n}$ generated by the vectors $C_i^{\rho}$ (see Sturmfels~
 \cite[Proposition 1.3]{sturmfelsEquations}).  Specifically, the orbit
 corresponding to a face $F$ of $P$ is the intersection of $V
 \subseteq \mathbb A^{nr}_{\K}$ with the subscheme $\{(b_i^\rho)\in
 \mathbb{A}^{nr}_{\K} : b_i^\rho\neq 0 \text{ for } C_i^{\rho} \in F,
 b_i^\rho = 0 \text{ otherwise}\}$, and the {\em distinguished point}
 of this orbit $(b_i^{\rho}) \in V$ satisfies $b_i^{\rho}=1$ for
 $C_i^{\rho} \in F$ and $b_i^{\rho}=0$ otherwise.  Note that the
 $T_C$-orbit of this distinguished point is the orbit associated with
 $F$. In particular, the standard torus $T_C \subseteq V$ corresponds
 to the full face $P$, and is the $T_C$-orbit of the distinguished
 point $(1,\dots,1) \in \mathbb A^{nr}_{\K}$.

For $\theta \in \Theta$, the torus orbits of $Y_{\theta}$ correspond
one-to-one to cones of $\mathcal N(P_{\theta})$ or, equivalently, to
faces of the polyhedron $P_{\theta}$.  Since $P_{\theta}=\pi_n(P \cap
\pi^{-1}(\theta))$, every face of $P_{\theta}$ is of the form
$F_{\theta}=\pi_n(\widetilde{F_{\theta}} \cap \pi^{-1}(\theta))$ where
$\widetilde{F_{\theta}}$ is the smallest face of $P$ containing the
preimage of $F_\theta$ under the map $\pi_n$.  The torus orbits of
$Y_{\theta}$ correspond one-to-one to the $\theta$-semistable
$T_C$-orbits in $V$, and distinguished points of $Y_{\theta}$
correspond one-to-one to $\theta$-semistable distinguished points of
$V$.

Recall from Corollary~\ref{cor:fan-of-normalization} that the support
of the fan $\mathcal N(P_{\theta})$ is $(\mathbb Q^n_{\geq 0})^\vee$.

\begin{definition}
  \label{defn:distinguishedgcon} For $\theta \in \Theta$ and
  $\mathbf{w} \in (\mathbb Q_{\geq 0}^n)^\vee$, the \emph{distinguished
  $\theta$-semistable representation} $b_{\theta,\mathbf{w}} =
  (b_i^\rho)$ is the distinguished point of $V$ corresponding to the
  cone of $\mathcal N(P_{\theta})$ containing $\mathbf{w}$ in its
  relative interior.  
\end{definition}


Computing these representations has been the key tool in understanding
the moduli space $\mathcal M_\theta$ (see \cite{Nakamura, Reid2,
  Crawthesis}).  No reasonable algorithm was known to compute these
representations for $\mathcal M_\theta \neq \ghilb$, and the algorithm
for \(\ghilb\) introduced by Nakamura required that one perform a
sequence of {\em $G$-igsaw transformations} to calculate a single
$G$-cluster.  These transformations are the \(\ghilb\) analogues of
{\em flips} for the toric Hilbert scheme introduced by
Maclagan--Thomas~\cite{MT1}.

 The next theorem presents an elementary method to compute any
 distinguished $\theta$-semistable quiver representation in one step.
 Note that the cone dual to \(P\) is \(P^\vee =
 \{(\mathbf{v},\mathbf{w}) \in (\Q^{r})^*\times (\Q^n)^* : w_i+v_\rho
 - v_{\rho\rho_i}\geq 0\}\).  Given $\mathbf{w} \in (\mathbb Q^n_{\geq
   0})^\vee$, consider the slice \(P^\vee_{\bf w} := \{\mathbf{v} \in
 (\Q^{r})^* : w_i+v_\rho - v_{\rho\rho_i}\geq 0\}\).

\begin{theorem}
 \label{thm:distinguishedcone}
 Fix $\theta \in \Theta$ and $\mathbf{w} \in (\mathbb Q^n_{\geq
   0})^\vee$.  Let $\mathbf{v} \in P^{\vee}_{\mathbf{w}}$ be any
 vector with $\theta \cdot \mathbf{v} \leq \theta \cdot \mathbf{v}'$
 for all $\mathbf{v}' \in P^{\vee}_{\mathbf{w}}$.  The distinguished
 $\theta$-semistable quiver representation $b_{\theta,{\bf
     w}}=(b_i^{\rho})$ has
 \begin{equation} \label{eqn:birho} b_i^\rho =
 \left\{\begin{array}{cl} 1 & \text{if } w_i+v_\rho - v_{\rho\rho_i} =
 0 \\ 0 & \text{if } w_i+v_\rho - v_{\rho\rho_i} > 0
 \end{array}\right. .
\end{equation} 

\end{theorem} 

\begin{proof} 
  Write $\face_{\bf w}(P_\theta)$ for the face of $P_\theta$ where
  ${\bf w}$ is minimized, and $F:= \widetilde{\face}_{\bf
    w}(P_\theta)$ for the smallest face of $P$ containing the preimage
  of $\face_{\bf w}(P_\theta)$ under the map $\pi_n$.  The distinguished point
  $b_{\theta,{\bf w}}\in V$ satisfies $b_i^\rho = 1$ if $C_i^\rho\in
  F$ and $b_i^\rho = 0$ otherwise. To restate this condition in terms
  of weight vectors, note that ${\bf v}\in P^\vee_{\bf w}$ implies
  $(\mathbf{v},\mathbf{w})\in P^\vee$, so $w_i+v_\rho - v_{\rho\rho_i}
  \geq 0$ for all \(1\leq i\leq n\) and \(\rho\in G^*\). Furthermore,
  $(\mathbf{v},\mathbf{w})\) lies in the face $\mathcal{N}_P(F)$ of
  $P^\vee$ if and only if $w_i+v_\rho - v_{\rho\rho_i} = 0$ for
  $C_i^\rho\in F$ and $w_i+v_\rho - v_{\rho\rho_i} > 0$ otherwise.  In
  particular, the distinguished point $b_{\theta,{\bf w}}\in V$
  satisfies the conditions of (\ref{eqn:birho}) if and only if
  $(\mathbf{v},\mathbf{w})\) lies in the face $\mathcal{N}_P(F)$.

We complete the proof by showing that for $\mathbf{v}$ satisfying the
hypothesis of the theorem, the vector $(\mathbf{v},\mathbf{w})$ lies
in the cone $\mathcal{N}_P(F)$.  The hypothesis states that ${\bf v}$
lies in $\face_{\theta}(P^\vee_{\bf w})$, the face of $P^\vee_{\bf w}$
where $\theta$ is minimized. Thus,
$(\mathbf{v},\mathbf{w})$ lies in the smallest face of $P^\vee$
containing $\face_{\theta}(P^\vee_{\bf w})$.  The fact that
$\mathcal{N}_P(F)$ is the smallest face of $P^\vee$ containing
$\face_{\theta}(P^\vee_{\bf w})$ is the content of
Craw--Maclagan~\cite[Proposition~2.7]{CrawMaclagan}.  Thus
$(\mathbf{v},\mathbf{w})\in \mathcal{N}_P(F)$.
\end{proof}

 \begin{remark}
   To compute the quiver representation $b\in V$ corresponding to a
   point $[b]\in Y_\theta$ that is not distinguished, first compute
   the distinguished point $b_{\theta, {\bf w}}$ in the same
   $T_C$-orbit as $b\in V$, and then let $T_C$ act on the coordinates
   of $b_{\theta, {\bf w}}$.
 \end{remark}

 \begin{example}
 \label{ex:wzero}
   For ${\bf w} = {\bf 0}\in (\mathbb Q^n_{\geq 0})^\vee$, we obtain the
   inequalities $v_\rho \geq v_{\rho \rho_i}$ for all $\rho \in G^*$,
   $1 \leq i \leq n$.  Since every arrow of the McKay quiver lies in
   some directed cycle, these inequalities must be equalities, so the
   quiver representation $b_{\theta, {\bf 0}}$ satisfies $b_i^\rho =
   1$ for all $1\leq i\leq n$ and $\rho\in G^*$.
 \end{example}

 \begin{remark}
   When ${\bf w}\in (\Q^n_{\geq 0})^\vee$ is a point of the lattice
   \(M^\vee\), the inequalities defining $P^\vee_{\bf w}$ form the
   \emph{reductor condition} of Logvinenko~\cite[Equation
   (4.8)]{Logvinenko}.  Thus, 
   Algorithm~\ref{alg:distinguished} enables Logvinenko to verify
   whether the $G$-constellations arising from his reductor sets are
   \(\theta\)-stable for any given $\theta$.
 \end{remark}

 Theorem~\ref{thm:distinguishedcone} gives an explicit algorithm to
 compute the distinguished McKay quiver representation $b_{\theta,
   {\bf w}}$ for $(\theta, {\bf w}) \in \Theta \times (\mathbb
 Q^n_{\geq 0})^{\vee}$, which we now present.

\begin{algorithm}
  \label{alg:distinguished} To compute the distinguished
  representation $b_{\theta,
   {\bf w}}$.
 \end{algorithm}

\noindent{\bf Input:} $(\theta, {\bf w}) \in \Theta \times (\mathbb
    Q^n_{\geq 0})^{\vee}$ and the matrix $C$.

 \medskip
\begin{enumerate} 
\item Compute the polyhedron $P^\vee_{\bf w} = \{ {\bf v'} \in \mathbb Q^r \,:\,
  ({\bf v'},{\bf w})C
    \geq {\bf 0} \}$.
\item Compute an optimal solution ${\bf v}$ of the linear program 
$$ \textup{minimize} \{ \theta \cdot {\bf v'} \,:\, {\bf v'} \in
P^\vee_{\bf w} \}.$$
\item The distinguished quiver representation $b_{\theta,{\bf
      w}}=(b_i^{\rho})$ has coordinates $$b_i^\rho =
  \left\{\begin{array}{cl} 1 & \text{if } w_i+v_\rho -
      v_{\rho\rho_i} =  0 \\
      0 & \text{if } w_i+v_\rho - v_{\rho\rho_i} > 0
 \end{array}\right..$$
\end{enumerate}

\begin{proof}[Proof of Correctness]
 This is immediate from Theorem~\ref{thm:distinguishedcone}.
\end{proof}



 \begin{example}
 \label{ex:hard}
 For the group action of type \(\frac{1}{11}(1,2,8)\) considered in
 Example~\ref{ex:128}, the given three-dimensional representation
 decomposes as $\rho_1\oplus\rho_2\oplus\rho_{8}$.
 We now compute a pair of
 \(\theta\)-stable distinguished quiver representations for the
 parameter $\theta = (1,1,1,1,-7,-9,1,1,1,8,1)$ considered in
 Example~\ref{ex:128}.
 
 The vector ${\bf w} = (10,7,6)$ lies in the relative interior of the
 cone generated by $(2,4,5)$, $(7,3,1)$ and $(1,0,0)$ corresponding to
 vertices 1, 4 and 6 in Figure~\ref{fig:128}.  The vector ${\bf v} =
 (-8,-10,-1, -3,6,4, -9, 0, -2, -15, -6)$ is an optimal solution to
 the linear program in Step (2) from
 Algorithm~\ref{alg:distinguished}.  Calculating the vector $({\bf
   v},{\bf w}) C = (12, 0, 0, 1, 0, 11, 12, 0, 11, 1, 0, 11, 12, 22,
 22, 23, 11, 11, 1, 0, 0, 12, 22, 0, 23, 11, 0, 1, 0, \\ 0, 12, 11,
 0)$, shows that the distinguished representation is $$
 b_{\theta,{\bf
     w}} = (0, 1, 1, 0, 1, 0, 0, 1, 0, 0, 1, 0, 0, 0, 0, 0, 0, 0, 0,
 1, 1, 0, 0, 1, 0, 0, 1, 0, 1, 1, 0, 0, 1), $$
 where the coordinates
 of $b_{\theta,{\bf w}} = (b_i^\rho)$ are indexed exactly as the
 columns of the matrix $C$.  Figure~\ref{fig:subquiver} shows only the
 arrows $a_i^\rho$ for which $b_i^\rho \neq 0$.

\begin{figure}[!ht]
\input{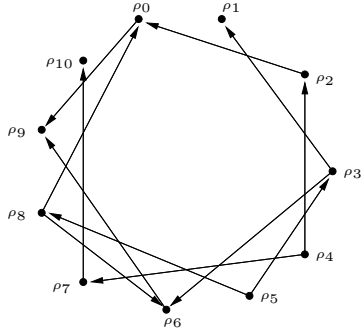}
  \caption{The subquiver consisting of arrows $a_i^\rho$ for which
    $b_i^\rho \neq 0$}
  \label{fig:subquiver}
\end{figure}

\end{example}

\begin{example}
  For the group action from Example~\ref{ex:hard}, consider the vector
  ${\bf w} = (8,3,1)$ in the relative interior of the two
  dimensional cone of $\mathcal{N}(P_{\theta})$ generated by the
  vectors $(1,0,0)$ and $(7,3,1)$ corresponding to vertices 6 and 4 in
  Figure~\ref{fig:128}.  In this case, the optimal solution is ${\bf
    v} = (-5, -9, -2, -6, 1,-3, -7, 0, -4, -8,-1)$.  The product
  $({\bf v},{\bf w}) C$ has 18 entries equal to zero, and $$
  b_{\theta,{\bf w}} = (0, 1, 1, 0, 1, 1, 0, 1, 1, 0, 1, 1, 0, 0, 0,
  0, 1, 1, 0, 1, 1, 0, 0, 1, 0, 1, 1, 0,1, 1, 0, 0,1).$$
\end{example}


\begin{thebibliography}{10}

\bibitem{BezrukavnikovKaledin}
R.~Bezrukavnikov and D.~Kaledin.
\newblock Mc{K}ay equivalence for symplectic resolutions of quotient
  singularities.
\newblock {\em Tr. Mat. Inst. Steklova}, 246:20--42, 2004.

\bibitem{Bollobas}
B.~Bollob{\'a}s.
\newblock {\em Modern graph theory}, volume 184 of {\em Graduate Texts in
  Mathematics}.
\newblock Springer-Verlag, New York, 1998.

\bibitem{BKR}
T.~Bridgeland, A.~King, and M.~Reid.
\newblock The {M}c{K}ay correspondence as an equivalence of derived categories.
\newblock {\em J. Amer. Math. Soc.}, 14(3):535--554 (electronic), 2001.

\bibitem{porta}
T.~Christof.
\newblock Porta, a software system to compute convex hulls.
\newblock Available by anonymous ftp from {\tt
  http://www.iwr.uni-heidelberg.de/groups/comopt/software/PORTA/}.

\bibitem{Crawthesis}
A.~Craw.
\newblock {\em The {M}c{K}ay correspondence and representations of the
  {M}c{K}ay quiver}.
\newblock PhD thesis, University of Warwick, (2001).

\bibitem{Craw}
A.~Craw.
\newblock An introduction to motivic integration.
\newblock In {\em Strings and geometry}, volume~3 of {\em Clay Math. Proc.},
  pages 203--225. Amer. Math. Soc., Providence, RI, 2004.

\bibitem{CrawIshii}
A.~Craw and A.~Ishii.
\newblock Flops of \({G}\)-\emph{{H}ilb} and equivalences of derived categories
  by variation of {GIT} quotient.
\newblock {\em Duke Math. J.}, 124(2):259--307, 2004.

\bibitem{CrawMaclagan}
A.~Craw and D.~Maclagan.
\newblock {Fiber fans and toric quotients}, (2005).
\newblock To appear in \emph{Discrete Comput. Geom.}, arXiv: math.AG/0510128.

\bibitem{CMT2}
A.~Craw, D.~Maclagan, and R.~R. Thomas.
\newblock {Moduli of McKay quiver representations II: Gr\"{o}bner basis
  techniques}, (2005).
\newblock preprint.

\bibitem{Dolgachev}
I.~Dolgachev.
\newblock {\em Lectures on Invariant Theory}, volume 296 of {\em London Math.
  Soc. Lecture Note Series}.
\newblock Cambridge University Press, Cambridge, 2003.

\bibitem{DGM}
M.~Douglas, B.~Greene, and D.~Morrison.
\newblock Orbifold resolution by {D}-branes.
\newblock {\em Nuclear Phys. B}, 506(1-2):84--106, 1997.

\bibitem{Fulton}
W.~Fulton.
\newblock {\em Introduction to toric varieties}, volume 131 of {\em Annals of
  Mathematics Studies}.
\newblock Princeton University Press, Princeton, NJ, 1993.

\bibitem{M2}
D.~Grayson and M.~Stillman.
\newblock Macaulay 2, a software system for research in algebraic geometry.
\newblock Available from {\tt http://www.math.uiuc.edu/Macaulay2/}.

\bibitem{Haiman}
M.~Haiman.
\newblock Hilbert schemes, polygraphs and the {M}acdonald positivity
  conjecture.
\newblock {\em J. Amer. Math. Soc.}, 14(4):941--1006 (electronic), 2001.

\bibitem{HaimanSturmfels}
M.~Haiman and B.~Sturmfels.
\newblock Multigraded {H}ilbert schemes.
\newblock {\em J. Algebraic Geom.}, 13(4):725--769, 2004.

\bibitem{HostenSturmfels}
S.~Ho{\c{s}}ten and B.~Sturmfels.
\newblock G{RIN}: an implementation of {G}r\"obner bases for integer
  programming.
\newblock In {\em Integer programming and combinatorial optimization
  (Copenhagen, 1995)}, volume 920 of {\em Lecture Notes in Comput. Sci.}, pages
  267--276. Springer, Berlin, 1995.

\bibitem{Ishii1}
A.~Ishii.
\newblock On the {M}c{K}ay correspondence for a finite small subgroup of {${\rm
  GL}(2,\mathbb C)$}.
\newblock {\em J. Reine Angew. Math.}, 549:221--233, 2002.

\bibitem{ItoNakajima}
Y.~Ito and H.~Nakajima.
\newblock Mc{K}ay correspondence and {H}ilbert schemes in dimension three.
\newblock {\em Topology}, \textbf{39}(6):1155--1191, 2000.

\bibitem{ItoNakamura}
Y.~Ito and I.~Nakamura.
\newblock Hilbert schemes and simple singularities.
\newblock In {\em New trends in algebraic geometry (Warwick, 1996)}, volume 264
  of {\em London Math. Soc. Lecture Note Ser.}, pages 151--233. Cambridge Univ.
  Press, Cambridge, 1999.

\bibitem{Kronheimer}
P.~Kronheimer.
\newblock The construction of {ALE} spaces as hyper-{K}\"ahler quotients.
\newblock {\em J. Differential Geom.}, 29(3):665--683, 1989.

\bibitem{Logvinenko}
T.~Logvinenko.
\newblock Families of $\rm{G}$-constellations over resolutions of quotient
  singularities.
\newblock arXiv:math.AG/0305194, (2003).

\bibitem{MT1}
D.~Maclagan and R.~R. Thomas.
\newblock Combinatorics of the toric {H}ilbert scheme.
\newblock {\em Discrete Comput. Geom.}, 27(2):249--272, 2002.

\bibitem{Mumford}
D.~Mumford, J.~Fogarty, and F.~Kirwan.
\newblock {\em Geometric invariant theory}, volume~34 of {\em Ergebnisse der
  Mathematik und ihrer Grenzgebiete (2)}.
\newblock Springer-Verlag, Berlin, 1994.

\bibitem{Nakamura}
I.~Nakamura.
\newblock Hilbert schemes of abelian group orbits.
\newblock {\em J. Algebraic Geom.}, 10(4):757--779, 2001.

\bibitem{PeevaStillman}
I.~Peeva and M.~Stillman.
\newblock Toric {H}ilbert schemes.
\newblock {\em Duke Math. J.}, 111(3):419--449, 2002.

\bibitem{Reid2}
M.~Reid.
\newblock Mc{K}ay correspondence.
\newblock In {\em Proc.\ of algebraic geometry symposium ({K}inosaki, {N}ov
  1996), {T}. {K}atsura (ed.)}, pages 14--41, (1997).

\bibitem{Reid3}
M.~Reid.
\newblock La correspondance de {M}c{K}ay.
\newblock {\em Ast\'erisque}, (276):53--72, 2002.
\newblock S\'eminaire Bourbaki, Vol.\ 1999/2000.

\bibitem{SI2}
A.~Sardo-Infirri.
\newblock Resolutions of orbifold singularities and the transportation problem
  on the {M}c{K}ay quiver.
\newblock arXiv: math.AG/9610005, (1996).

\bibitem{Schrijver}
A.~Schrijver.
\newblock {\em Theory of linear and integer programming}.
\newblock The Wiley-Interscience Series in Discrete Mathematics. John Wiley \&
  Sons Ltd., Chichester, 1986.

\bibitem{GBCP}
B.~Sturmfels.
\newblock {\em Gr\"obner bases and convex polytopes}, volume~8 of {\em
  University Lecture Series}.
\newblock American Mathematical Society, Providence, RI, 1996.

\bibitem{sturmfelsEquations}
B.~Sturmfels.
\newblock Equations defining toric varieties.
\newblock In {\em Algebraic geometry---Santa Cruz 1995}, volume~62 of {\em
  Proc. Sympos. Pure Math.}, pages 437--449. A.M.S., Providence, RI, 1997.

\bibitem{Thompson}
H.~Thompson.
\newblock {Fan is to monoid as scheme is to ring: a generalization of the
  notion of a fan}, (2003).
\newblock arXiv:math.AG/0306221.

\end{thebibliography}
 \end{document}